\documentclass[amsthm]{elsart}
\usepackage{amsmath,amssymb,amsfonts,amscd,graphicx,epigraph,texdraw}
\usepackage[all]{xy}

\newcommand{\proj}{\mathbb P}
\newcommand{\aff}{\mathbb A}

\newcommand{\kay}{{\mathbf k}}
\newcommand{\Kay}{{\mathbb K}}
\newcommand{\Eff}{{\mathbb F}}

\newcommand{\com}{\mathbb{C}}
\newcommand{\zee}{{\mathbb{Z}}}

\newcommand{\real}{{\mathbb{R}}}
\newcommand{\rat}{\mathbb{Q}}
\newcommand{\oh}{{\mathcal{O}}}
\newcommand{\ca}{{\mathcal{A}}}

\newcommand{\cf}{{\mathcal{F}}}

\newcommand{\cc}{{\mathcal{C}}}

\newcommand{\cz}{{\mathcal{Z}}}

\newcommand{\cw}{{\mathcal{W}}}
\newcommand{\ccv}{{\mathcal{V}}}

\newcommand{\cy}{{\mathcal{Y}}}

\newcommand{\cd}{{\mathcal{D}}}
\newcommand{\cb}{{\mathcal{B}}}

\newcommand{\ccr}{{\mathcal{R}}}

\newcommand{\ccq}{{\mathcal{Q}}}

\newcommand{\cu}{{\mathcal{U}}}

\newcommand{\cx}{{\mathcal{X}}}

\newcommand{\mfm}{{\mathfrak{m}}}

\def\<{\ensuremath{\langle}}
\def\>{\ensuremath{\rangle}}

\newcommand{\Chow}{{\operatorname{Chow}}}
\newcommand{\Hilb}{{\operatorname{Hilb}}}
\newcommand{\Univ}{{\operatorname{Univ}}}

\newcommand{\HI}{{\operatorname{HI}}}
\newcommand{\CI}{{\operatorname{CI}}}
\newcommand{\rdim}{{\operatorname{rdim}}}
\newcommand{\supp}{{\operatorname{supp}}}
\newcommand{\rank}{{\operatorname{rank}}}
\newcommand{\Vol}{{\operatorname{Vol}}}

\newcommand{\Conv}{{\operatorname{Conv}}}

\newcommand{\Trop}{{\operatorname{Trop}}}

\newcommand{\Span}{{\operatorname{Span}}}
\newcommand{\Spec}{{\operatorname{Spec}}}
\newcommand{\Hom}{{\operatorname{Hom}}}

\def\ol#1{{\overline{#1}}}
\def\init#1#2{{{\operatorname{in}}_#1(#2)}}

\newtheorem{property}[thm]{Property}

\begin{document}
\setcounter{secnumdepth}{0}
\setcounter{tocdepth}{ 3 }

\begin{frontmatter}
\title{A Tropical Toolkit}
\author{Eric Katz}
\address{Department of Mathematics\\ University of Texas\\eekatz@math.utexas.edu}
\begin{abstract}
{\small We give an introduction to Tropical Geometry and prove some
results in Tropical Intersection Theory. The first part of this
paper is an introduction to tropical geometry aimed at researchers
in Algebraic Geometry from the point of view of degenerations of
varieties using projective not-necessarily-normal toric varieties.
The second part is a foundational account of tropical intersection
theory with proofs of some new theorems relating it to classical
intersection theory.}
\end{abstract}

    \epigraph{O monumento \'{e} bem moderno.}{Caetano Veloso \cite{Veloso}}
\end{frontmatter}  

\numberwithin{thm}{section}

  \tableofcontents
 
    \section{Introduction}

  Tropical Geometry is an exciting new field of mathematics
    arising out of computer science.  In the mathematical realm,
    it has been studied by Mikhalkin
    \cite{Mik06}, Speyer \cite{Speyer}, the Sturmfels school \cite{Sturm03}, Itenberg, Kharlamov, and Shustin \cite{IKS},
    Gathmann and Markwig \cite{GMCH}, and Nishinou and Siebert \cite{NS}
    among many others.  It has found applications in the enumeration of curves \cite{Mik03}, low-dimensional topology \cite{BSLL}, algebraic dynamics \cite{EKL}, and the study of compactifications \cite{HKT,Tevelev}.

    This paper is an introduction to tropical geometry from the
    point of view of degenerations of subvarieties of a toric
    variety.  In this respect, its approach is close to that of
    the Sturmfels school.

    In the first part of the paper, we use
    not-necessarily-normal projective toric varieties to introduce
    standard notions such as degenerations, the
    Gr\"{o}bner and fiber fans, and tropical varieties.  In the
    second part of the paper, we give a foundational account of
    tropical intersection theory.  We define tropical
    intersection numbers, and show that tropical
    intersection theory computes classical intersection numbers under certain hypotheses,
    use tropical intersection theory to get data on deformation of
    subvarieties, and associate a tropical cycle to subvarieties.  The two parts can be read independently.  

    We will express tropical geometry in the language of projective not-necessarily
    normal toric schemes over a valuation ring (see \cite{GKZ}, Chapter 5 for such
    toric varieties over fields).  These
    toric schemes give toric degenerations.  There are other
    constructions of toric degenerations analogous to different constructions of toric varieties.  
    Analogous to the fan construction as in \cite{FultonToric} is the approach of Speyer \cite{Speyer}.  See also the paper of Nishinou and Siebert \cite{NS}.  In \cite{Uniformizing}, Speyer introduced a construction of toric degenerations paired with a map to projective space.      The construction we use here has
    the advantage of being very immediate at the expense of some loss of
    generality by mandating projectivity and the loss of computability versus more constructive methods.
    
    We have chosen in this paper to approach the material from the
    point of view of algebraic geometry and had to neglect the very beautiful combinatorial
    nature of this theory.  We would like to suggest that the reader takes a look at \cite{Sturm03} for a more
    down-to-earth introduction to tropical geometry.  We also point out a number of
    references that are more combinatorial in nature and which relate to our approach.  There is
    the wonderful book of Gelfand, Kapranov, and Zelevinsky
    \cite{GKZ} which gives a combinatorial description of the
    secondary polytope among many other beautiful results, the paper of
    Billera and Sturmfels on fiber polytopes \cite{FiberPoly} (see
    also the lovely book of Ziegler \cite{Polytopes}), the book of
    Sturmfels on Convex Polytopes and Gr\"{o}bner Bases \cite{GBaCP}
    as well as the papers \cite{KSZ,GrobnerToric}.

    We should mention that since this paper first appeared in preprint form, there has emerged a synthetic approach to tropical intersection theory.  The intersection theory of tropical fans was established by Gathmann, Kerber, and Markwig in \cite{GKM} and was extended to general tropical varieties in $\real^n$ by Allermann and Rau in \cite{AR}.
  
    Many of the results from the first part of this paper are rephrased from Speyer's
    dissertation \cite{Speyer} and the general outlook is implicit
    in the work of Tevelev \cite{Tevelev} which introduced the interplay between toric
    degenerations and tropical compactifications.  Please see
    \cite{DFS} for an explanation of the relationship between such
    work.  We hope this piece will be helpful to other researchers.

    We would like to thank Bernd Sturmfels for suggesting the
    connection between tropical cycles and Minkowski weights
    and Hannah Markwig, David Speyer, Fr\'{e}d\'{e}ric Bihan, and Sam Payne for helpful comments and corrections.

    \section{Conventions}

    Let $\ccr$ be a ring with a valuation contained in a subgroup $G$ of $(\real,+)$,
    $$v:\ccr\setminus\{0\}\rightarrow G\subseteq \real.$$
    Let $\Kay$ denote the field of fractions of $\ccr$, $\mfm$ the maximal ideal $v^{-1}((0,\infty))$, and
    $\kay=\ccr/\mfm$.

    There are two examples that will be most important:
    \begin{enumerate}
    \item[(1)] $\Kay=\com\{\{t\}\}=\bigcup_M \com((t^{\frac{1}{M}}))$, the field of formal Puiseux series,
    $v:\Kay\rightarrow\rat$, the order map and $\kay=\com$.

    \item[(2)] $\Kay=\com((t^{\frac{1}{M}}))$, the field of formal Laurent
    series in $t^{\frac{1}{M}}$, $v:\Kay\rightarrow\zee\frac{1}{M}$, and $\kay=\com$.
    \end{enumerate}

    Note that the first choice of $\ccr$ has the disadvantage of
    not being Noetherian.  This is not much of a hindrance because
    any variety defined over $\Kay$ in the first case can be defined over
    $\Kay$ in the second case for some $M$.  This will be enough in practice.

    In either case, given $x\in\Kay$, we may speak of the {\em leading term} of $x$.  This is the non-zero complex coefficient of the lowest power of $t$ occurring in the power-series expansion  of $x$.

    In either of these cases we have an inclusion
    $\kay\hookrightarrow\ccr$ such that the composition
    \[\kay\hookrightarrow\ccr\rightarrow\ccr/\mfm=\kay\]
    is the identity.

    Also, for every $u\in G$, we have an element $t^u\in\ccr$
    so that $v(t^u)=u$.  These elements have the property
    that
    \[t^{u_1}t^{u_2}=t^{u_1+u_2}.\]
   The choice of a map $u\mapsto t^u$ as a section of $v$ is perhaps unnatural.  In \cite{PayneFibers}, Payne introduced a formalism of tilted rings which avoids the need for a section.
   
    For an $n$-tuple, $w=(w_1,\dots,w_n)\in G^n$, we may write
    $t^w$ for $(t^{w_1},\dots,t^{w_n})\in(\Kay^*)^n$.  Similarly,
    we may write $v:(\Kay^*)^n\rightarrow G^n$ for the product of
    valuations.

    For $g=(g_1,\dots,g_n)\in(\Kay^*)^n$, $\chi=(\chi_1,\dots,\chi_n)\in\zee^n$, we write $g^\chi$ for
    $g_1^{\chi_1}\dots g_n^{\chi_n}\in\Kay^*$.

    \section{Polyhedral geometry}
    
    Here we review some notions from polyhedral geometry.  Please see \cite{Polytopes} for more details.  
    
    Let $\ca\subset\real^n$ be a set of points.  Let $P=\Conv(\ca)$ be their convex hull.  For $v\in(\real^n)^\vee$, the face $P_v$ of $P$ is the set of points $x\in P$ that minimize the function $\<x,v\>$.  Let $\Gamma_v=\ca\cap P_v$.  The cone 
    \[C_{\Gamma_v}=\left\{w\in(\real^n)^\vee\ \Big|\begin{array}{cccl}
    \ \<\chi_i,w\>&=&\<\chi_j,w\>& \text{ for } \chi_i,\chi_j\in\Gamma_v\\
    \ \<\chi,w\>&<&\<\chi',w\>& \text{ for } \chi\in\Gamma_v,\chi'\notin\Gamma_v
    \end{array} \right\}\]
    is the normal cone to the face $P_v$.  Observe that $v$ is in the relative interior of $C_{\Gamma_v}$.  
    The correspondence between $P_v$ and $C_{\Gamma_v}$ is inclusion reversing.  The $C_\Gamma$'s form a fan, $N(P)$, called the (inward) normal fan of $P$.  
    
    Two polytopes are said to be normally equivalent if they have the same normal fan.  
    
    A polyhedron in $\real^n$ is said to be {\em integral} with respect to a full-rank lattice $\Lambda\subset (\real^n)^\vee$ if it is the intersection of half-spaces defined by equations of the form $\{x|\<x,w\>>a\}$ for $w\in\Lambda$, $a\in\real$.  We will usually not note the lattice when it is understood.
    
    \begin{defn}
    A {\em polyhedral complex} in $\real^n$ is a finite collection $\cc$ of polyhedra in $\real^n$ that contains the faces of any one of its members, and such that any non-empty intersection of two of its members is a common face.
    \end{defn}
    
    A polyhedral complex is said to be {\em integral} if all of its members are integral polyhedra.  The support $|\cc|$ of a polyhedral complex $\cc$ is the set-wise union of its polyhedra.  We say that a polyhedral complex $\cc$ is supported on a polyhedral complex $\cd$ if $|\cc|\subseteq|\cd|$.
    
    \begin{defn}
    Given two integral polyhedral complexes, $\cc$,$\cd$ in
    $\real^n$, we say $\cc$ is a refinement of $\cd$ if every polyhedron in $\cd$ is a union of polyhedra in $\cc$.
     \end{defn}

    It is well-know that for convex polytopes $P$ and $Q$ with
    normal fans $N(P)$, $N(Q)$, $N(P)$ is a refinement of
    $N(Q)$ if and only $\lambda Q$ is a Minkowski summand of $P$ for some $\lambda\in\real_{>0}$.
    See Proposition 1.2 of
    \cite{BFS}.
    
    Given a polyhedron $P$ in a complex $\cc$, we may construct a fan $\cf$ called the {\em star} of $P$.  Pick a point $w$ in the relative interior of $P$  Let $\cd$ be the set of all polyhedra in $\cc$ containing $P$ as a face.  For every $Q\in\cd$, let $C_Q$ be the cone
    \[C_Q=\{v\in\real^n| w+\epsilon v\in Q \text{ for some } \epsilon>0\}.\]
    These $C_Q$'s give a fan $\cf$.  If $P$ is a maximal polyhedron in $\cc$, then its star is its affine span.
   Please note that this usage of star is non-standard.
    
    \begin{defn} Given $n$ polytopes, $P_1,\dots,P_n\subset \real^n$, their {\em mixed volume} is the coefficient of $\lambda_1\lambda_2\dots\lambda_n$ of
    $\Vol(\lambda_1P_1+\dots+\lambda_nP_n)$
    which is a homogeneous polynomial of degree $n$ in $\lambda_1,\dots,\lambda_n$.
    \end{defn}
    
    \section{Toric Schemes}

    \subsection*{Toric Schemes over $\Spec\ \ccr$}
     
    We take the point of view of \cite{Smi97} and use the language of toric
    schemes over $\Spec\ \ccr$.  We
    use the not necessarily normal projective toric varieties of
    \cite{GKZ}.

    For $T=(\Kay^*)^n$ a $\Kay$-torus, let $T^{\wedge}=\Hom(T,\Kay^*)$ be the
    character lattice and $T^\vee=\Hom(\Kay^*,T)$ be the
    one-parameter subgroup lattice.  Let
    $T^\wedge_\real=\real\otimes T^\wedge$, $T^\vee_\real=\real\otimes
    T^\vee$, and $T^\vee_G=G\otimes T^\vee$.

    A homomorphism of tori $T\rightarrow U$ induces homomorphisms
    $T^\vee\rightarrow U^\vee$ and $U^\wedge\rightarrow T^\wedge$.

    \begin{defn}
    Let $T=(\Kay^*)^n\hookrightarrow (\Kay^*)^{N+1}/(\Kay^*)\hookrightarrow
    PGl_{N+1}(\Kay)$ be a composition of homomorphisms of groups
    where $(\Kay^*)^{N+1}/(\Kay^*)$ denotes the quotient by the diagonal
    subgroup and the last homomorphism is the diagonal inclusion.  For
    $y\in\proj^N_{\Kay}$, let $T_y$ denote the stabilizer of $y$ in $T$.  The toric
    variety associated to $(T,y)$ is the closure
    $$Y=\ol{(T/T_y)y}.$$
    $Y$ lies in the fiber over the generic point in $\proj^N_{\ccr}\rightarrow\Spec\ \ccr$.  Let
    the toric scheme $\cy$ be the closure of $Y$ in $\proj^N_\ccr$, and let $Y_0=\cy\times_{\Spec\
    \ccr}\Spec\ \kay$ be the special fiber.
    \end{defn}

    \begin{defn} If $y\in\proj^n_{\kay}\subset\proj^n_{\Kay}$ for $\kay\subset
    \Kay$ then the toric scheme is said to be {\em defined over
    $\kay$}.  Alternatively, it's obtained by base-change from a
    toric variety defined over $\kay$ by the map $\Spec\ \Kay\rightarrow\Spec\ \kay$
    induced by the inclusion.
    \end{defn}

    \begin{exmp} \rm
    \label{firstsegre}
    Let $T=(\Kay^*)^2\rightarrow (\Kay^*)^4/(\Kay^*)$ be the
    inclusion given by
    $$(x_1,x_2)\mapsto (1,x_1,x_2,x_1x_2).$$
    If $y=[1:1:1:1]\in\proj^3_\Kay$ then
    $$T\cdot y=\{[1:x_1:x_2:x_1x_2]\ |\ x_1,x_2\in\Kay^*\}.$$
    The closure of the above is $\proj^1\times\proj^1$ under the
    Segre embedding.  This is defined over $\kay$.
    \end{exmp}

    \begin{defn} There is a natural map from $(\Kay^*)^n$ to $Y$ given by
    $$\xymatrix{
    (\Kay^*)^n\ar[r]      &Y\\
    g      \ar@{|->}[r]&g \cdot y
    }$$
    The image of the map is called {\em the big open torus}.
    If the map is an open immersion, we say our toric variety is
    {\em immersive}.
    \end{defn}

    Now, we explain a method of defining toric schemes.
    Let $\ca=\{\chi_1,\dots,\chi_{N+1}\}\subset T^\wedge=\zee^n$ be a finite set.  Let $a:\ca\mapsto G$
    be a function called a {\em height function}.
    Let ${\mathbf{y}}=(y_1,\dots,y_{N+1})\in(\Kay^*)^{N+1}$ be
    an element satisfying
    \[v(y_i)=a(\chi_i).\]
    The choice of $\ca$ induces a homomorphism of groups
    \[\begin{array}{rll}
    T=(\Kay^*)^n&\rightarrow& (\Kay^*)^{N+1}\\
    g=(g_1,\dots,g_n)&\mapsto&(g^{\chi_1},\dots,g^{\chi_{N+1}}).
    \end{array}\]
We may consider the map as a homomorphism
    $T\rightarrow(\Kay^*)^{N+1}/(\Kay^*)$ where the quotient is
    by the diagonal subgroup.  Therefore, if $\mathbf{y}\in(\Kay^*)^{N+1}$,
    \[g\cdot y=(g^{\chi_1}y_1,\dots,g^{\chi_{N+1}}y_{N+1}).\]

    One may ask how the toric variety depends on the choice of $y$. Let ${\mathbf{y,y}}'\in (\Kay^*)^{N+1}$ satisfy
    \[v(y_i)=v(y'_i)=a(\chi_i).\]
    Then $\mathbf{y,y'}$ are related by multiplication by an element $g\in(\Kay^*)^{N+1}$ with $v(g)=0$.  This element lifts to an element of $({\mathbb G}_m)_\ccr^{N+1}$.  Therefore, the two choices of $\cy_{\ca,a}$ are related by an action of the diagonal torus in $\proj^N_\ccr$.  As a consequence, the special fibers are related by an action of the diagonal torus in $\proj^N_\kay$.  

    Let $\cy_{\ca,a}$ be the toric scheme associated to $T$ and $y$.  Note that
    if the integral affine span of $\ca$ is $\zee^n$ then $Y_{\ca,a}$ is immersive.

    It is a theorem that the normalization of $Y$ is the toric variety
    associated to the normal fan of the polytope $\Conv(\ca)$.  See \cite{Cox} for details.

    \begin{defn} The {\em induced subdivision} of
    $\Conv(\ca)$ is given as follows.  Let the {\em upper hull} of
    $a$ be
    $$\operatorname{UH}=\Conv(\{(\chi,b)|\chi\in\ca,\ b\geq
    a(\chi)\}).$$
    The faces of $\operatorname{UH}$ project down to give a subdivision of
    $\Conv(\ca)$.
    \end{defn}
      
     $\Conv(\ca)$ is called the {\em weight polytope} of $Y$ while the induced subdivision is called the {\em weight subdivision} of $\cy$.

    \begin{exmp} \rm
    \label{latticesquare}
    Let $\ca=\{(0,0),(1,0),(0,1),(1,1)\}$ be the
    vertices of a lattice square.  Let $a$ be given by
    $$a(0,0)=0,\ a(1,0)=0,\ a(0,1)=0,\ a(1,1)=1.$$
    Choose ${\mathbf y}=(1,1,1,t^1)$. 
    This induces the inclusion
    $T\hookrightarrow(\Kay^*)^4/(\Kay^*)$ given by
    $$(x_1,x_2)\mapsto
    (x_1^0x_2^0,x_1^1x_2^0,x_1^0x_2^1,x_1^1x_1^1)=(1,x_1,x_2,x_1x_2)$$
    as in Example \ref{firstsegre}.
    Therefore $\cy$ is the closure of the image of
    $$(x_1,x_2)\mapsto [1:x_1:x_2:tx_1x_2].$$
    The fiber over $\Spec\ \Kay$ is isomorphic to the closure of
    $$(x_1,x_2)\mapsto [1:x_1:x_2:x_1x_2],$$
    which is $\proj_\Kay^1\times\proj_\Kay^1$ under the Segre embedding.

    The special fiber can be seen as follows: taking the limit of $(x_1,x_2)$ as
    $t\mapsto 0$, we get $[1:x_1:x_2:0]$ which is $\proj^2$; taking the
    limit of $(t^{-1}x_1,t^{-1}x_2)$ as $t\mapsto 0$, we get
    $[0:x_1:x_2:x_1x_2]$ which is another $\proj^2$.  One sees that the
    special fiber is two copies of $\proj^2$ joined along
    $\proj^1$.  We will show that this case is indicative of a general phenomenon in Lemma \ref{torusstructure}.
    \end{exmp}

    \begin{figure}
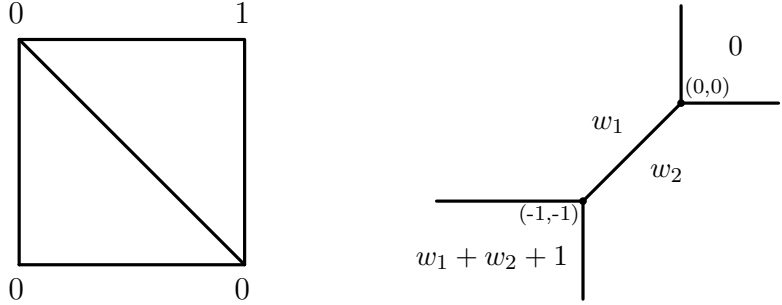


\begin{center}
   \begin{texdraw}
       \drawdim cm  \relunitscale 1.5
       \linewd 0.03
        \move(0 0) \lvec (2 0)
        \lvec (2 2)
        \lvec( 0 2)
        \lvec(0 0)
        \move (2 0) \lvec(0 2)
        \htext(-.1 -.3) {$0$}
        \htext(1.9 -.3) {$0$}
        \htext(-.1 2.15) {$0$}
        \htext(1.9 2.15) {$1$}
     \end{texdraw}     
     \hspace{2cm}
     \begin{texdraw}
       \drawdim cm  \relunitscale 1.3
       \linewd 0.03
       \move(0 0) \lvec(-1.5 0)
       \move(0 0) \lvec(0 -1)
      \move (0 0) \fcir f:0 r:0.04
      \htext(-.65 -.25){\scriptsize (-1,-1)}  
      \move(0 0) \lvec(1 1)
      \lvec(1 2)
      \move(1 1) \lvec(2 1)
      \move (1 1) \fcir f:0 r:0.04
      \htext(1.05 1.05) {\scriptsize (0,0)}
   
     \htext(-1.7 -.7) {\footnotesize $w_1+w_2+1$}
     \htext(.7 .2) {\footnotesize $w_2$}
     \htext(.1  .7) {\footnotesize $w_1$}
     \htext(1.5 1.5){\footnotesize $0$}
   \end{texdraw}
\end{center}

    \caption{The subdivision and its dual complex}
    \label{square}

    \end{figure}

    \subsection*{Recovering the Weight Subdivision}

    There is a way of working backwards from $(T,y)$ to $\ca$ and
    a subdivision of $\Conv(\ca)$.

    \begin{defn} Let $V$ be a $\Kay$-vector space.  A $\kay$-weight decomposition is a
    vector space isomorphism defined over $\kay\subset\Kay$
    $$V\cong\bigoplus_{\chi\in\zee^n} V_\chi$$
    where $H$ acts on $V_\chi$ with character $\chi$.
    \end{defn}

    \begin{lem} Any $\Kay$-vector space $V$ on which $H$ acts
    linearly has a $\kay$-weight decomposition.
    \end{lem}

    \begin{proof} See \cite{Borel}, Propositions 8.4 and 8.11.\end{proof}

    Lift $y\in\proj^N_{\Kay}$ to ${\bf y}\in\Kay^{N+1}$
    Write ${\bf y}=\sum_\chi {\bf y_\chi}$.  Let $\ca=\{\chi\in
    \zee^n|y_\chi\neq 0\}$.  Then $\Conv(\ca)$ is called the
    weight polytope of $Y$.
    If $\dim V_\chi=1$, set
    $a_\chi=v({\bf v_\chi})$.
    Otherwise, write ${\bf v_\chi=v_1+\dots+v_n}$ where ${\bf v_i}$ are
    vectors in a one-dimensional subspace on which $H$ acts, and
    set $a_\chi=\min(v({\bf v_i)})$.
    Take the subdivision of $\Conv(\ca)$ induced by $a_\chi$ which is independent of the lift $y$.

    \subsection*{Dual Complex}

    Consider the pairing
    \[T_\real^\wedge\otimes T^\vee_\real\rightarrow\real\]
    and the piecewise-linear function
    \[F:T^\vee_\real\rightarrow\real\]
    defined by
    \[F(w)=\min_{\chi\in\ca} (\<\chi,w\>+a_\chi).\]
    The domains of linearity of $F$ give a polyhedral complex structure on $T^\vee_\real$.
    For $\Gamma\subset\ca$, let
     \[C_{\Gamma}=\left\{w\in(\real^n)^\vee\ \Big|\begin{array}{cccl}
    \ \<\chi_i,w\>+a_{\chi_i}&=&\<\chi_j,w\>+a_{\chi_j}& \text{ for } \chi_i,\chi_j\in\Gamma\\
    \ \<\chi,w\>+a_{\chi}&<&\<\chi',w\>+a_{\chi'}& \text{ for } \chi\in\Gamma,\chi'\notin\Gamma
    \end{array} \right\}\]
    If $C_\Gamma$ is not empty, then $\Gamma$ are points of $\ca$ in a face of the weight subdivision.  The $C_\Gamma$'s fit together to form an integral polyhedral complex,  the  {\em dual complex} which is dual to the weight subdivision.  Note that if $a_\chi=0$ for $\chi\in\ca$, the weight subdivision becomes the weight polytope and the dual subdivision becomes the normal fan.
    \begin{exmp} \rm Figure \ref{square} shows the weight subdivision and
    dual complex for Example \ref{latticesquare}.
    Here,
    \[\begin{array}{rcl}
    F(w)&=&\min_{\chi\in\ca}(\<\chi,w\>+a_\chi)\\
    &=&\min(0,w_1,w_2,w_1+w_2+1).
   \end{array}\]
   The values of $F$ on the dual complex are noted in the figure.    
    \end{exmp}

    \subsection*{One-parameter Families of Points}

    Let us review the notion of specialization.  For $y\in \proj_\Kay^N$, we may take $\ol{y}\in \proj_\ccr^N$, considered as a scheme over $\Spec\ \ccr$.  The {\em specialization} of $y$ is
    \[\hat{y}=\ol{y}\times_{\Spec\ \ccr}\Spec\ \kay\in\proj^N_\kay.\]
    We can compute the specialization by hand.  Lift $y$ to ${\mathbf y}\in\Kay^{N+1}\setminus\{0\}$ such that $\min(v(y_i))=0$.   If $y_i\neq 0$, write $y_i=c_1t^{b_i}+\dots$ where the ellipsis denotes higher order terms.
    Let
\[S=\{i|b_i=0\}\]

     Then $\mathbf {\hat{y}}$ satisfies

     \[\hat{y_i}=\begin{cases}c_i &\text{if } i\in S\\
     0 &\text{else.}\end{cases}\]

    \begin{defn}   Let $\cy$ be a toric scheme over $\ccr$.
    Let $y$ be a point in $Y$.  Given $g\in (\Kay^*)^n$, the family associated to $(g,y)$ is
    the scheme over $\Spec\ \ccr$ given by the closure of $g\cdot y$.
    \end{defn}

    \begin{defn} The limit of $(g,y)$ is the point in $Y_0$
    given by
    $$\ol{g\cdot y}\times_{\Spec\ \ccr} \Spec\ \kay$$
    \end{defn}

    Now, observe that 
    \[v((g\cdot y)_i)=<\chi_i,v(g)>+a_{\chi_i}\]
    where $v(y_i)=a_{\chi_i}$.
    Therefore, when we base-change to $\Spec\ \kay$, the only components of $g\cdot y_i$ that stay non-zero are the ones on which $<\chi_i,v(g)>+a_{\chi_i}$ is minimized.
    Consequently, if $v(g)\in C_\Gamma$ for a cell $\Gamma$ of the weight subdivision, and $\hat{y}$ is the limit of $(g,y)$, then $\hat{y}_i\neq 0$ if and only if $\chi_i\in\Gamma$.

    \subsection*{One-parameter Families of Subschemes}
     
   We will also consider degenerations of subschemes $X$ of $Y$.

    \begin{defn}
    Let $w\in G^n$ and  $g=t^w$.
    Consider the subscheme of $\cy$ given by $\ol{g\cdot X}$, the closure of
    $g\cdot X$.  Define the {\em initial degeneration} of $X$ to be the subscheme of $Y_0$ given by
    $$\init{w}{X}=\ol{g\cdot X}\times_{\Spec\ \ccr}\Spec\ \kay$$
    \end{defn}

    \begin{exmp} \rm This definition specializes to the usual
    definition of the initial form of a polynomial.  Let
    \[f=x_1^2x_2+7x_1x_2x_3+4x_3^3\in\Kay[x_1,x_2,x_3],\]
    and let $w=(3,4).$ Let $X=V(f)\subset Y=\proj^2_\Kay$.
    Then $t^wV(f)$, a subvariety of $\proj^2_\Kay$ is $V(h)$ for
    \[\begin{array}{rcl}
    h&=&(t^{-3}x_1)^2(t^{-4}x_2)+7(t^{-3}x_1)(t^{-4}x_2)(x_3)+4(x_3)^3\\
     &=&t^{-10}x_1^2x_2+7t^{-7}x_1x_2x_3+4x_3^3\\
     &=&t^{-10}(x_1^2x_2+7t^3x_1x_2x_3+4t^{10}x_3^3).
    \end{array}\]
    Therefore,
    $$\init{w}{V(f)}=\ol{t^w\cdot V(f)}\times_{\Spec\ \ccr}\Spec\
    \kay$$
    is cut out by
    $$\init{w}{f}=x_1^2x_2.$$
    \end{exmp}

    Now that if $X=x$ is a point then, $\init{w}{X}=\ol{t^w\cdot
    x}\times_{\Spec\ \ccr}\Spec\ \kay$.

    Every point of $\init{w}{X}$ occurs as a limit of the form
    $\ol{g\cdot x} \times_{\Spec\ \ccr} \Spec\ \kay$ for $x\in X$.  This is the content of the {\em tropical lifting lemma}.   This lemma was first announced without proof in \cite{Systems}.  A proposed proof was given in \cite{TropGras} but has been found to be incomplete.  A proof using affinoid algebras was given by Draisma in \cite{Draisma}.  Jensen, Markwig, and Markwig provided an algorithm that finds a tropical lift in \cite{JMM}.  This algorithm uses some ideas from our proof and their paper is recommended as an exposition of our proof in terms of commutative algebra.  In \cite{PayneFibers}, by applying a projection argument to reduce to the hypersurface case, Payne gave a stronger version of tropical lifting that works over more general fields. 
    
    We first review the concept of relative dimension from Chapter 20 of \cite{Fulton}.  
    
    \begin{defn} \label{relatdim} Let $p:Z\rightarrow S$ be a scheme over a regular base scheme $S$.  For $V$, a closed integral subscheme of $Z$,  let $T=\ol{p(V)}$.  
   The {\em relative dimension}
    of $V$ is
    \[\rdim\ccv=\operatorname{tr.deg.}(R(V)/R(T))-\operatorname{codim}(T,S).\]
    \end{defn}
         
    We will apply this definition for $T=\Spec\ \com[[t^{\frac{1}{M}}]]$.  Note that a point in the special fiber is of relative dimension $-1$.

    \begin{lem} [Tropical Lifting Lemma]
    \label{limitpoints} Let $\Kay=\com\{\{t\}\}$.
    If $\tilde{x}\in\init{w}{X}$ then there exists
    $x\in X$ with
    $$\init{w}{x}\equiv \ol{t^w\cdot x}\times_{\Spec\ \ccr}\Spec\ \kay=\tilde{x}.$$
    \end{lem}

    \begin{proof}
    We treat $X$ as a subscheme of $\proj^N_\Kay$.
    If $\dim X=0$, then the support of $t^wX$ is a union of closed $\Kay$-points.  One such point specializes to $\tilde{x}$.  The corresponding component has initial
    deformation supported on $\tilde{x}$ and gives the desired point in $X$. Therefore, we may suppose $\dim X=n>0$.

    Pick $M$ sufficiently large so that $X$ is defined over
    $\Eff=\com((t^{1/M}))$.  Let $\ccq=\com[[t^{1/M}]]$.

    By replacing $X$ by $t^{w}X$ we may
    suppose $w=0$. Let $\ol{X}$ be the closure of $X$ in $\cy$.
    Note that $\ol{X}$ is flat over $\Spec\ \ccq$.

    Let $W_0$ be a codimension $n$ subvariety of $Y_0\subset \proj_k^N$ such that
    $W_0$ intersects
    $$X_0=\ol{X}\times_{\Spec\ \ccq} \Spec\ \kay$$
    in a
    $0$-dimensional subscheme containing $\tilde{x}$.  Extend
    $W_0$ to a flat integral scheme $\cw\rightarrow\Spec\ \ccq$ so that
    $\cw\times_{\Spec\ \ccq}{\Spec\ \kay}=W_0$ (for example, we may set $\cw=W_0\times_{\Spec\ \kay}\Spec\ \ccq$).  Then,
    $\ol{X}\times_{\cy}\cw$ is a scheme, all of whose components have non-negative relative dimension 
    over $\Spec\ \ccq$.
    The following equality holds for underlying sets
    $$(\ol{X}\times_{\cy}\cw)\times_{\Spec\ \ccq}{\Spec\
    \kay}=X_0\times_{Y_0} W_0.$$
    Since the scheme on the right is $0$-dimensional, there are no
    components of $\ol{X}\times_{\cy}\cw$ contained in the
    special fiber.  Therefore, the induced reduced structure on $\ol{X}\times_{\cy}\cw$ is
    flat, has relative dimension $0$ and has a component of its limit
    supported on
    $\tilde{x}$.  Let $W=\cw\times_{\Spec\ \ccq}\Spec\ \Eff$.  By
    uniqueness of flat limits, the closure of the induced reduced
    structure on $X\times_Y W$ in $\cy$ is the induced reduced structure on
    $\ol{X}\times_{\cy}\cw$.

    Therefore, we
    may apply the $0$-dimensional case to the induced reduced structure on $X\times_Y
    W$.
    \end{proof}
       
   We will find the following corollary useful.
      
    \begin{cor} Under the hypotheses of the previous lemma
    and the equality of underlying sets $X=\ol{X\cap(\Kay^*)^n}$,
    we may suppose $x\in X\cap (\Kay^*)^n$.
    \label{openlimitpoints}
    \end{cor}

    \begin{proof}
    Produce $x\in X$ as above.  If $x\in X\cap (\Kay^*)^n$ then we
    are done.  Otherwise, there is a morphism
    \[f:\Spec\ \Kay[[s]]\rightarrow X\]
    so that the generic point is sent to $X\cap (\Kay^*)^n$ while
    the closed point is sent to $x$.  This morphism is defined
    over some $\com((t^{\frac{1}{M}}))$ and can be given as a base-change from
    \[f:\Spec\ \com[[t^{\frac{1}{M}}]][[s]]\rightarrow X\]
    where we view $X$ as defined over $\com[[t^{\frac{1}{M}}]]$.
    Therefore, we may extend the morphism to $f:\Spec\ \com[[t^{\frac{1}{M}}]][[s^{\frac{1}{M}}]]\rightarrow X$.
    Consider the diagonal morphism
    $$i:\Spec\ \com[[u^\frac{1}{M}]]\rightarrow\Spec\
    \com[[t^{\frac{1}{M}}]][[s^{\frac{1}{M}}]]$$
    induced by
    $$t^{\frac{1}{M}}\mapsto u^{\frac{1}{M}},\ s^{\frac{1}{M}}\mapsto
    u^{\frac{1}{M}}.$$
    By restricting the composition $f\circ i$ to the generic point, $\Spec\
    \com((u^{\frac{1}{M}}))$, we find the desired $\Kay$-point.
    \end{proof}

    \subsection*{Structure of $\cy_{\ca,a}$}

    $\cy_{\ca,a}$ has well-understood fibers over the generic and special
    point.

    \begin{defn} For $\Gamma$, a face of the weight polytope,
    let $Y^0(\Gamma)\subset Y$ be the set of all points $y\in Y\subseteq \proj_{\Kay}^N$
    so that their lifts ${\mathbf y}\in(\Kay)^{N+1}\setminus\{0\}$
    satisfy
    \[{\bf y_i}\neq 0 \text{ if and only if } \chi_i\in\Gamma\]
    \end{defn}

    \begin{defn} For $\Gamma$, a cell of the weight
    subdivision, let $Y_0^0(\Gamma)\subset Y_0\subset \proj_{\kay}^N$
    be the set of all points $y\in Y_0\subseteq \proj_{\kay}^N$
    so that their lifts ${\mathbf y}\in(\kay)^{N+1}\setminus\{0\}$
    satisfy
    \[{\bf y_i}\neq0 \text{ if and only if } \chi_i\in\Gamma.\]
     \end{defn}

    \begin{prop} \label{torusstructure}
    \begin{enumerate}
    \item[]
    \item[(1)] $Y=\cy_{\ca,a}\times_{\Spec\ \ccr} \Spec\ \Kay$ is the toric
    variety associated to $\ca$.  The non-empty faces of the weight polytope are in
    inclusion-preserving bijective correspondence with its torus
    orbits given by $\Gamma\mapsto Y^0(\Gamma)$.

    \item[(2)] The scheme $Y_0=\cy_{\ca,a}\times_{\Spec\ \ccr} \Spec\ \kay$
    is supported on the union of toric varieties associated to the top-dimensional cells of the
    weight subdivision such that the
    non-empty cells of the weight subdivision are in
    inclusion-preserving bijective correspondence with its torus
    orbits given by $\Gamma\mapsto Y^0_0(\Gamma)$.
    \end{enumerate}
    \end{prop}

    \begin{proof}
    (1) is Proposition 1.9 of Chapter 5 of
    \cite{GKZ}.
    We give the proof of (2) which is directly analogous.  Elements of
    $\cy_{\ca,a}\times_{\Spec\ \ccr} \Spec\ \kay$ are of the form
    \[\ol{g\cdot y}\times_{\Spec\ \ccr}\Spec\ \kay\] by
    Lemma \ref{openlimitpoints}.  If $v(g)\in C_\Gamma$, the cell of the dual complex corresponding to $\Gamma$, then the
    limit $\ol{g\cdot y}\times_{\Spec\ \ccr}\Spec\ \kay$ is in the
    orbit $Y^0_0(\Gamma)$.

    Similarly if $w\in C_\Gamma$, by varying $g$ with $v(g)=w$, we may make 
    $\ol{g\cdot y}\times_{\Spec\ \ccr}\Spec\ \kay$ be any point of $Y_0^0(\Gamma)$.
    \end{proof}

    Part (2) of the above lemma is simply not true at the level of
    scheme structure.  As a counterexample, take $\ca=\{0,1,2\}$,
    $a(0)=0,\ a(1)=1,\ a(2)=0$.  Then $Y_0$ is a double-line in
    $\proj^2$.  The corresponding subdivision is the single cell
    $[0,2]$ whose toric variety is the reduced-induced structure
    on $Y_0$.  The construction of toric degenerations by fans as in
    \cite{Speyer} is better behaved in this respect.

    In the case of Example \ref{latticesquare}, we see that $Y_0$ consists of two $\proj^2$'s, five $\proj^1$'s and four fixed-points.

    It is instructive to phrase the above theorem in the language
    of the dual complex.  Given two elements $g,g'\in (\Kay^*)^n$ with
    $v(g)=v(g')$,
    the limits of $(g,y)$ and
    $(g,y')$ are related by the action of an element of
    $(\kay^*)^n$ and so lie in the same open torus orbit.  Therefore, we
    may define an equivalence relation on $G^n$.  Two elements
    $w,w'\in T_G^\vee$ are equivalent, written $w\sim_y w'$ if for $g,g'\in G$ satisfying $w=v(g)$ and $w'=v(g')$, 
    the limits of $(g,y)$ and $(g',y)$ lie in the same open torus
    orbit.

    \begin{prop} \label{dualequiv}
    $w\sim_y w'$ if and only if $w$ and
    $w'$ lie in the same cell in the dual complex associated to the
    toric scheme $\cy_{\ca,a}$.
    \end{prop}
    

    \subsection*{Invariant Limits}

   The open orbits $Y^0(\Gamma)$ and $Y_0^0(\Gamma)$ are fixed point-wise by sub-tori in $T$.
   
   \begin{lem} \label{fixedtori} Let $\Gamma$ be a face of the weight polytope (resp. cell of the weight subdivision).  Let $w\in C_\Gamma$ and $H\subset T$ be the sub-torus with $H_\real^\vee=\Span(C_\Gamma-w)$.  Let $z\in Y^0(\Gamma)$ (resp. $Y^0_0(\Gamma)$). Then the maximal sub-torus fixing $z$ is $H$.
  \end{lem}
  
  \begin{proof} We give the proof for $Y^0_0(\Gamma)$.  The proof for $Y_0(\Gamma)$ is similar.
  
  Let $z\in Y^0_0(\Gamma)$.  Lift $z$ to $\mathbf{z}\in \kay^{N+1}\setminus\{0\}$.
  Every $g\in H$ satisfies $g^{\chi_i}=g^{\chi_j}$ for $\chi_i,\chi_j\in\Gamma$.   Let $g'=g^\chi\in\kay^*$ for $\chi\in\Gamma$.  Since $z_i\neq 0$ if and only if $\chi_i\in\Gamma$,
  \[g\cdot {\mathbf{z}}=g'\mathbf{z}\]
  which is another lift of $z$,
  
  If $u\in T^\vee\setminus H^\vee$, there exists $\chi_i,\chi_j\in\Gamma$ such that 
  \[\<\chi_i,u\>\neq \<\chi_j,u\>.\]
  it follows is that $z$ is not fixed by the one-parameter subgroup corresponding to $u$.
  \end{proof}

    We may rephrase the above lemma.
    \begin{lem}
    \label{InvariantTorus}
    Suppose $g\in T$ satisfies $v(g)\in C_\Gamma$.  Then the limit of 
    $(g,y)$ in $Y_0$ is invariant under the torus $H$ given by $H_\real^\vee=\Span(C_\Gamma-w)$.
    \end{lem}

    \begin{proof}
    The limit of $(g,y)$ lies in $Y_0^0(\Gamma)$.
    \end{proof}

    Suppose $v(g)$ lies in $C_\Gamma$, the cell of the dual complex dual to a cell $\Gamma$
    in the weight subdivision.
    We may make use of the map $\Spec\ \ccr\rightarrow\Spec\
    \kay$ to base-change the limit
    $$\ol{g\cdot y}\times_{\Spec\ \ccr} \Spec\ \kay$$
    to
    $$\hat{y}=(\ol{g\cdot y}\times_{\Spec\ \ccr} \Spec
    \ \kay)\times_{\Spec\ \kay} \Spec\ \ccr.$$
    This just means that we should consider a limit point's coordinates as points in $\Kay$ rather than in
    $\kay$ and take its closure.

    \begin{lem} The weight polytope of the toric scheme
    $\widehat{\cy}=\ol{(\Kay^*)^n\cdot \hat{y}}$ is $\Conv(\Gamma)$.
    \label{ConstantCoeff}
    \end{lem}

    \begin{proof}
    Lift $\hat{y}$ to $\mathbf{\hat{y}}\in \Kay^{N+1}\setminus\{0\}$.  The weights with which $T$ acts on $\mathbf{\hat{y}}$ are $\chi\in\Gamma$.  Therefore the weight polytope in $\Conv(\Gamma)$.
    \end{proof}

    The dual complex of $\widehat{\cy}$ is the normal fan of $\Conv(\Gamma)$.  One may also that the normal fan of $\widehat{\cy}$ is the star of $C_\Gamma$, the cell of the dual complex dual to $\Gamma$.

    \begin{lem}\label{smalld}  Let $\hat{y}=\init{w}{X}$.  For $u\in T^\vee_G$, 
    \[\init{u}{\hat{y}}=\init{{w+\epsilon u}}{y}\]
    for sufficiently small $\epsilon>0$.
     \end{lem}

    \begin{proof}
     Let $w\in C_\Gamma$ for  $\Gamma$, a cell of the weight subdivision.  Then the weight polytope of $\hat{y}$ is $\Conv(\Gamma)$.  Therefore, $u$ is in a cone of the normal fan of $\Gamma$ dual to some face $\Gamma'\subseteq\Gamma$.  It follows that the coordinates of $\init{u}{\hat{y}}$ in $\proj^N$ are non-zero only for $\chi_i\in\Gamma'$ and in that case are equal to the leading terms of the coordinate of $t^wy$.   Now, $C_{\Gamma}$ is a face of $C_{\Gamma'}$ and we may pick small $\epsilon>0$ such that $w+\epsilon u\in C_{\Gamma'}$.  Therefore, $\init{{w+\epsilon u}}{y}=\init{u}{\hat{y}}$.
    \end{proof}

    \subsection*{Naturality of Dual Complexes}

    \begin{lem}
    \label{Refinement}
    Given a proper surjective $(\Kay^*)^n$-equivariant morphism of
    $n$-dimensional toric schemes, $f:\cx\rightarrow\cy$ then the dual complex of
    $\cx$ is a refinement of that of $\cy$.  The normal fan to the weight polytope of $X$ is a refinement of that of $Y$.
    \end{lem}

    \begin{proof} Let $x\in \proj_{\Kay}^N$ so that
    $\cx=\ol{T\cdot x}$ and
    $\cy=\ol{T\cdot f(x)}$ for (possibly different) diagonal actions
    of $T$ on $\proj_{\Kay}^N, \proj_{\Kay}^{N'}$.   
    
    Now let $C_\Gamma$ be a $k$-dimensional cell in the dual complex of $\cx$.  We must show that $f^\vee(C_\Gamma)$ is in the relative interior of a cell in the dual complex of $\cy$ of dimension at least $k$.  If $g\in T$ satisfies $v(g)\in C_\Gamma$, then the limit, $\hat{x}$ of $(g,x)$ is invariant under the $k$-dimensional torus $H$ with $H_\real^\vee=\Span(C_\Gamma-v(g))$.  Since $f$ is equivariant, $f(\hat{x})$ is the limit of $(g,f(x))$.  Furthermore if $v(g)\in C_{\Gamma'}$, a cell in the dual complex of $\cy$ then $f(\hat{x})$ is invariant under an $l$-dimensional torus $H'$ with ${H'}_\real^\vee=\Span(C_{\Gamma'}-v(g))$.  Since $f(\hat{x})$ is also invariant under $H$, then $l>k$.  
    
    To prove the statement for the weight polytope, we may set $\cx=X\times\Spec\ \Kay[[s]]$, $\cy=Y\times\Spec\ \Kay[[s]]$ where $s$ is an algebraic indeterminate.  Consider the valuation $v:\Kay[[s]]\rightarrow\zee$ given by $v(s)=1$, $v(\Kay^*)=0$.  Then the weight subdivision of $\cx$ and $\cy$ are exactly the weight polytopes of $X$ and $Y$ and the same argument applies.
    \end{proof}

    \subsection*{Equivariant Inclusions}

    In this section we consider a projection of integral polytopes
    $p:P\rightarrow Q$ where $P=\Conv(\ca)$.

    \begin{defn}
    Given a finite set $\ca\subseteq \zee^n$ and a function
    $$a:\ca\rightarrow\real,$$
    a projection $p:\zee^n\rightarrow\zee^m$,
    let $\cb=p(\ca)$ and define {\em the image height function}
    $$b:\cb\rightarrow\real$$
    by
    $$b(\psi)=\min(\{a(\chi)|\chi\in p^{-1}(\psi)\}).$$
    The associated subdivision is the {\em image subdivision}.
    \end{defn}

    Note that the image subdivision is dependent on the height
    function not just on the original subdivision.  Weight polytopes and weight subdivisions are contravariant.

    \begin{lem}
    \label{EquivariantInclusion}
    Let $i:T\hookrightarrow
    U$ be an injective homomorphism of tori, so
    $$\ol{T\cdot v}\hookrightarrow \ol{U\cdot v}\hookrightarrow \proj^n$$
    is a chain of equivariant inclusions.  Then the induced
    projection
    $$i^\wedge:U^\wedge\rightarrow T^\wedge$$
    takes the weight polytope and the weight subdivision of $\ol{U\cdot v}$ to those of
    $\ol{T\cdot v}$.
    \end{lem}

    \begin{proof} The proof is straightforward.
    \end{proof}

    \section{Degenerations}

    \subsection*{Moduli Spaces}
    
    Tropical geometry is, in a certain sense, a method of
    parameterizing degenerations of subvarieties of a toric
    variety.  There are two useful moduli spaces for
    parameterizing degenerations, the Chow variety and the Hilbert
    scheme.   Points in these moduli spaces correspond to cycles or to subschemes.  This is useful because limits of points in the moduli space correspond to limits of cycles and subschemes.  This allows us to apply the machinery developed in the previous section to limits of subvarieties.
    
    Let $Y\subseteq \proj^N$ be a projective toric
    variety whose torus action extends to one on $\proj^N$.
    Recall that $k$-dimensional algebraic cycles of $Y$ are finite formal sums of $k$-dimensional subvarieties of $Y$ with integer coefficients.    Consider a subvariety $X\subset Y$, with degree $d$ in $Y$ and
    Hilbert polynomial $P$.  There are two moduli spaces that one can construct, $\Chow_d(Y)$ and $\Hilb_P(Y)$ that each have a point corresponding to $X$.  Points in $\Chow_d(Y)$ correspond to cycles in $Y$ of degree $d$.  We denote the point (called the {\em Chow form}) in $\Chow_d(Y)$ corresponding to $X$ by $R_X$.  $\Chow_d(Y)$ is constructed as a closed subscheme of $\Chow_d(\proj^N)$ which is a projective scheme.  Points in $\Hilb_P(Y)$ correspond to closed subschemes of $Y$ with Hilbert polynomial $P$.  The point $[X]$ in $\Hilb_P(Y)$ corresponding to $X$ is called the {\em Hilbert point}.
    Similarly, $\Hilb_P(Y)$ is a closed subscheme of $\Hilb_P(\proj^N)$ which is projective.
    
    See \cite{Kollar} for an
    in-depth construction of both varieties.  See also \cite{GKZ}
    for a discussion of the Chow variety.  We will
    break from the usage in \cite{Kollar} and use $\Chow$ to
    denote the un-normalized Chow variety which is there called
    $\Chow'$.  Note that the Hilbert scheme can
    be constructed over an arbitrary Noetherian scheme $S$ while
    there are restrictions on the base-scheme of the Chow variety.

    Let us review some useful properties of the Chow varieties and Hilbert schemes.
    
    \begin{property} The torus action on $Y$ induces a group action on
    $\Chow_d$ and $\Hilb_P$ which extends to an action on the ambient projective space.
    \end{property}
    
    Because the torus $T$ acts on $Y$, for $g\in T$, $g\cdot X$ is a subvariety of $Y$ of degree $d$ and Hilbert polynomial $P$.  Therefore, $R_{g\cdot X}\in\Chow_d(Y)$ and $[g\cdot X]\in\Hilb_P(Y)$.  This induces $T$-actions on $\Chow_d(Y)$ and $\Hilb_P(Y)$ given by
    \[\begin{array}{rclcrcl}
    T\times \Chow_d(Y)&\rightarrow&\Chow_d(Y),&\ \ \ &T\times\Hilb_P(Y)&\rightarrow&\Hilb_P(Y)\\
    (g,R_X)&\mapsto&R_{g\cdot X},&\ \ \ &(g,[X])&\mapsto&[g\cdot X]
    \end{array}\]
    
    \begin{property} There is a natural equivariant morphism
    $FC:\Hilb_P\rightarrow \Chow_d$ (see 5.4 of \cite{MFK}) called the
    {\em fundamental class map} that takes a scheme to its underlying cycle.
     \end{property}
     
     A subscheme $X$ of $Y$ has an underlying cycle.  Therefore, one may define a map 
     \[\begin{array}{rcl}
     FC:\Hilb_P&\rightarrow& \Chow_d\\
     \ [X]&\mapsto&R_X.
     \end{array}\]
This map is equivariant with respect to the above $T$-actions.
     
     \begin{property} The Hilbert scheme possesses a universal flat
    family $\Univ_P\rightarrow\Hilb_P$.
      \end{property}
      
     This universal family $\Univ_P$ is a subscheme of $Y\times\Hilb_P(Y)$.  The fiber over the Hilbert point $[X]$ is the subscheme $X$.  In particular, if $\Spec\ \Kay\rightarrow\Hilb_P(Y)$ is the $\Kay$-point $[X]$ then $\Univ_P\times_{\Hilb_P(Y)}\Spec\ \Kay=X$.
   
   The Chow variety does not usually have a universal flat family.  
   
   \begin{property} The Hilbert scheme is natural under base-change.  If  $Y\rightarrow S$ is projective then $\Hilb_P(Y/S)$ parameterizes $S$-subschemes of $Y$ with Hilbert polynomial $P$.  If $Z\rightarrow S$ is a morphism then
    \[\Hilb_P(Y\times_S Z/Z)=\Hilb_P(Y/S)\times_S Z.\]
    \end{property}
    
   The Chow variety does not have this property.
    
    The Hilbert scheme with its universal flat family and naturality properties is a much better behaved moduli space.  This makes it more useful for our purposes.    However, there are very beautiful combinatorial structures associated with the Chow variety.  See \cite{GKZ} for details.
   
  Now, we may use the Hilbert scheme to relate deformations of subschemes to limits of the form $(g,y)$.  Let $X$ be a subscheme of a toric variety $Y$.  Let $g\in T$ and $w=v(g)$.   By uniqueness of flat limits, the $\Spec\ \ccr$-point $\ol{g\cdot [X]}$ is the Hilbert point of $\ol{g\cdot X}$ in $\Hilb_P(\cy)$. 
  Therefore, the specialization of $g\cdot [X]$,
  \[\ol{g\cdot [X]}\times_{\Spec\ \ccr}\Spec\ \kay\in\Hilb_P(\cy)\times_{\Spec\ \ccr}\Spec\ \kay=\Hilb_P(Y_0)\]
 is the Hilbert point, $[\ol{g\cdot X}\times_{\Spec\ \ccr}\Spec\ \kay]$.    We may pull back the universal family by $\Spec\ \ccr\rightarrow\Hilb_P(\cy)$ to get a scheme $\cu$
    over $\Spec\ \ccr$.      
    Its special fiber is $\ol{g\cdot X}\times_{\Spec\ \ccr}\Spec\ \kay$.   If $g=t^w$, then the special fiber is the initial degeneration $\init{w}{X}$.

    \subsection*{Associated Toric Schemes}

    Let $Y$ be a toric scheme in $\proj_\Kay^N$ with a
    torus $T$.  Let $X$ be a subvariety of $Y$.  We may take the Hilbert point $[X]\in\Hilb_P(Y)$ or the
    Chow form $R_X\in\Chow_d(Y)$ and consider the two toric
    schemes, called the Hilbert and Chow images, respectively
    $$\HI=\ol{T/T_X\cdot[X]}\subseteq\Hilb_P(Y),\ \ \CI=\ol{T/T_X\cdot
    R_X}\subseteq\Chow_d(Y)$$
    where $T_X$ denotes the stabilizer of $[X]$ or $R_X$.

    \begin{defn} The subdivisions (in $(T/T_X)^\wedge\subseteq T^\wedge$) associated to the Hilbert
    and Chow images are called the {\em state subdivision} and the {\em secondary subdivision},
    respectively.  The dual polyhedral complexes (in $(T/T_X)^\vee$) are
    called the {\em Gr\"{o}bner complex} and the {\em Chow complex}.  In the case where
    $X$ and $Y$ are defined over $\kay$, these notions become the
    {\em state polytope, fiber polytope, the Gr\"{o}bner fan, and
    the fiber fan}, respectively.
    \end{defn}
    
      In the case where $X$ is also
    a toric subvariety in $Y$, the name fiber polytope is standard.  Otherwise our usage is somewhat
    non-standard.
  
    Now we may apply Proposition \ref{dualequiv} to the Gr\"{o}bner complex.
    
    \begin{prop}\label{dualequivgrob} Two points $w,w'$ lie in the same cell in the Gr\"{o}bner complex if and only if $\init{w}{X}$ and $\init{w'}{X}$ are related by a $T_\kay$-action.
    \end{prop}
    
    In the case where $X$ is defined over $\kay$, this proposition is close to the usual definition of the Gr\"{o}bner fan.   The usual definition, however, is a refinement
    of our definition.  This is because the initial ideals in the standard definition are sensitive to embedded primes associated to the irrelevant ideal.  Our definition is not.  The definition we give is based on that of  \cite{BayerMorrison}.    
  
   We may also apply Lemma \ref{InvariantTorus} to the Gr\"{o}bner complex.

    \begin{lem} If $w\in G^n$ is in the relative interior of a
    $k$-dimensional cell of the Gr\"{o}bner complex of $X$ then the closed
    subscheme $\init{w}{X}$ is invariant under a $k$-dimensional
    torus.
    \end{lem}

    \begin{proof}
    By Lemma \ref{InvariantTorus} the Hilbert point of
    $\init{w}{X}$ is invariant under a $k$-dimensional torus.
    Therefore, the closed subscheme $\init{w}{X}$ is invariant
    under the same torus.
    \end{proof}
 
    \begin{lem}\label{smalldef}
    For $u\in T^\vee_G$, 
    \[\init{u}{\init{w}{X}}=\init{{w+\epsilon u}}{X}\]
    for $\epsilon>0$ sufficiently small.
    \end{lem}
    
    \begin{proof}
    This is Lemma \ref{smalld} applied to the Hilbert point $[X]$.
    \end{proof}

    There is a natural projection
    $p:T_\real^\vee\rightarrow(T/T_X)_\real^\vee$.  We may abuse notation and
    use the term Gr\"{o}bner or Chow complex to also denote the
    appropriate complex's inverse image under $p$.

    \begin{exmp} \label{hypersurface} \rm
    Let $Y$ be a toric variety defined over $\kay$ given by a set of exponents $\ca\subset\zee^n$.  Let $X$ be a hypersurface defined in $Y$ by
    \[f(x)\equiv\sum_{\omega\in\ca} a_\omega x^\omega=0\]
    where $a_\omega\in\Kay$ and $x^\omega$ are coordinates on $Y\subset\proj^{|\ca|-1}$.  
    We may treat $[a_\omega]$ as coordinates on a projective space $(\proj^{|\ca|-1})^\vee$.  The torus $T$ acts on $(\proj^{|\ca|-1})^\vee$ by 
    \[\begin{array}{rcl}
    T\times (\proj^{|\ca|-1})^\vee&\rightarrow&(\proj^{|\ca|-1})^\vee\\
    (g,[a_\omega])&\mapsto &[g^{-\omega}a_\omega].
    \end{array}\]
    Then the equation
    \[\sum_{\omega\in\ca} a_\omega x^\omega=0\] 
    cuts out a universal hypersurface $\cu\subset Y\times(\proj^{|\ca|-1})^\vee$ over $(\proj^{|\ca|-1})^\vee$.  This universal family is flat and therefore defines a $T$-equivariant morphism  $(\proj^{|\ca|-1})^\vee\rightarrow\Hilb_P(Y)$.  The image of this morphism contains the Hilbert point of $X$.  Therefore, the Hilbert image, $T\cdot [X]$ is isomorphic to $Y$ but with the opposite torus action.  Therefore, the state polytope, which is the weight polytope of the Hilbert image is $-\Conv(\ca)$.  The Gr\"{o}bner fan is the normal fan $N(-\Conv(\ca))$.
 
For a down-to-earth exposition of this example, see Proposition 2.8 of \cite{GBaCP}.
    \end{exmp}

    \begin{exmp} \rm
    \label{ofapoint}
    Suppose that $Y$ is a toric variety defined over $\kay$.  Let
    $X$ be a reduced $\Kay$-point contained in an open torus orbit $Y^0(\Gamma)$.  The Hilbert scheme parameterizes reduced points in $Y$.  Therefore, the Hilbert image is $Y(\Gamma)$, the closure of $Y^0(\Gamma)$.  The weights on the Hilbert point of $X$ are $\chi\in\Gamma$ while the height function is  $a(\chi)=\<\chi,v(X)\>$.  It follows that the piecewise-linear function $F$ whose domains of linearity are the cells of the dual complex is
    \[F(w)=\min_{\chi\in\Gamma}\<\chi,w+v(X)\>.\]
    In particular
    if $x$ lies in the big open torus of $Y$ then the Gr\"{o}bner
    complex is just the normal fan of $\Conv(\Gamma)$ translated by $-v(X)$.

    Let us examine initial
    deformations if $X$ is a point in the big open torus in a toric variety $Y$.
    If $w=-v(X)$ then $t^wX$ has valuation $0$ and so $\init{w}{X}$
    is a point in the big open torus of $Y_0$.  Otherwise,
    $\init{w}{X}$ lies in some torus orbit.  In fact, if $w+v(X)\in C_\Gamma$ for a face $\Gamma$ of $Y$'s polytope then
    $\init{w}{X}$ is a point in $Y^0(\Gamma)$.  This is in agreement with the proof of Proposition \ref{torusstructure}.    
    \end{exmp}

    \begin{exmp} \rm
    Let $Y$ be a toric variety defined over $\kay$.
    Let $X$ be the scheme-theoretic image of a map $\Spec\ \kay[\epsilon]/
    \epsilon^2\rightarrow Y$.  We visualize $X$ as a point in $Y$ with a tangent vector anchored at it. Suppose the image lies in the big open
    torus and that the vector is chosen generically.
    Let us find the weight polytope of $\HI$.  By Proposition \ref{torusstructure}, it suffices to
    find the vertices which correspond to the torus-fixed points
    in $\HI$.  The torus-fixed points in $\HI$ are schemes $S$ consisting of a fixed point $p$ of $Y$ together with a projectivized tangent vector pointing along a $1$-dimensional orbit $E$ containing $p$.  By the genericity condition, all choice of $(p,E)$ with $p\in E$ are possible.  We must find the weights corresponding to these fixed points.
    
    Let us first work out the case where $Y=\proj^n$.
    If $\HI\subset \proj^N$ and ${\bf y}\in\kay^{N+1}\setminus\{0\}$ is a
    vector corresponding to a torus fixed point $Q$, then the
    vertex of the weight polytope of $\HI$ corresponds to the character of the
    action of $T=(\kay^*)^n$ on ${\bf y}$.  Because the
    embedding of $\HI$ is given by the composition of the embedding
    of the Hilbert scheme into a Grassmannian with the Pl\"{u}cker
    embedding into $\proj^N$, the action of $T$ on ${\bf y}$ is the same as
    the action of $T$ on
    $\wedge^{\text{top}}(\Gamma(\oh_Q(l)))$ where $l$ is a
    sufficiently large positive integer.   Now, a torus fixed-point of $\HI$ consists of a pair $(p,E)$.  Suppose $p$ is given by the point $X_i=\delta_{ir}$ in homogeneous coordinates.  Let $x_j=\frac{X_i}{X_r}$ be inhomogeneous coordinates on $X_r\neq 0$.  Then the fixed point $Q$ is given as the image of an affine morphism 
   \[\begin{array}{rcl} 
    \aff^n&\leftarrow&\Spec\ \kay[\epsilon]/\epsilon^2\\
    \kay[x_1,\dots,\hat{x}_r,\dots,x_{n+1}]&\rightarrow&\kay[\epsilon]/\epsilon^2\\
    x_i&\rightarrow&c\delta_{is}\epsilon
    \end{array}\]
    where $c\in\kay$ is some constant.  In other words, the tangent vector points along the $x_s$-axis.
    The vector space $\oh_Q(l)$ is spanned by two monomials, $X_r^l$ and $X_r^{l-1}X_s$.  They have characters $le_r$ and $(l-1)e_r+e_s$, respectively where $e_i$ are the standard unit basis vectors of $T^\wedge$.  Therefore, $\wedge^{\text{top}}(\Gamma(\oh_Q(l)))$ has character $(2l-2)e_r+(e_r+e_s)$.  Let $\Delta^{n-1}$ be the unit simplex in $T^\wedge$ and $\Gamma$ the convex hull of the mid-points of $2\Delta$.  Then the state polytope of $X$ which is the weight polytope of $\HI$ is  $(2l-2)\Delta+\Gamma$.  
    
 For a general toric variety $Y\subseteq\proj^n$, we note that the Hilbert scheme $\Hilb_P(Y)$ is constructed as a subscheme of $\Hilb_P(\proj^n)$.   Let $U$ be the torus of $\proj^n$, $T$ the torus of $Y$, $i:T\rightarrow U$ the homomorphism of tori, and $i^\vee:U^\vee\rightarrow T^\vee$ the induced projection.  If $Q$ is a $T$-fixed point of $\Hilb_P(Y)$, then $Q$ is a $U$-fixed point and the character of the corresponding vertex in $U^\vee$ pulls back by $i^\vee$ to the appropriate character in $T^\vee$.  Therefore, if $\Gamma=\Conv(\ca)$ is the polytope corresponding to $Y$ and $\Delta$ the convex hull of the mid-points of the edges of $2\Gamma$, the state polytope of $X$ is $(2l-2)\Gamma+\Delta$ by Lemma \ref{EquivariantInclusion}.    See \cite{OnnSturmfels} for a computation of the related case of the Gr\"{o}bner fan of generic point configurations in affine space. 
      
    The Chow image in this case is isomorphic to $Y$ as its points correspond to
    points of $Y$ with multiplicity $2$.  The fiber polytope
    is $\Gamma$.  Because the fiber polytope, $P$ is a Minkowski summand of the state polytope 
    $(2l-2)\Gamma+\Delta$, the Gr\"{o}bner fan is a
    refinement of the fiber fan.  This is an example of a general
    fact.
    \end{exmp}
   
    \begin{prop} The Gr\"{o}bner complex is a refinement of the
    fiber complex.
    \end{prop}

    \begin{proof}
    The fundamental class map $FC:\HI\rightarrow \CI$ satisfies
    the hypotheses of Lemma \ref{Refinement}.
    \end{proof}
    For a combinatorial commutative algebra proof of the above, see
    \cite{GrobnerToric}.
   
    \section{Tropical Varieties}

    \subsection*{Intersection of Sub-tori}

   Before we give the definition of tropical varieties,  we must digress to consider the intersection two sub-tori in $(\kay^*)^n$. Let
    \[H_1=(\kay^*)^{m_1},H_2=(\kay^*)^{m_2}\hookrightarrow T=(\kay^*)^n\]
    be two injective homomorphisms with $m_1+m_2=n$ such that images under the induced
    maps $H_i^\vee\rightarrow T^\vee$ are transversal.  Let
    $y_1,y_2\in(\kay^*)^n$.
    Let $V_i=H_i\cdot y_i$.  We compute the intersection
    of $V_1$ and $V_2$.

    The inclusions $H_1,H_2\hookrightarrow (\kay^*)^n$ correspond to
    surjections
    $T^\wedge\rightarrow H_i^\wedge$.
    Let $M_i$ be the kernel of the surjections.  We may also write
    $M_i$ as $H_i^\perp$.

    \begin{prop}
    \label{intersectsubtori}
    The number of intersection points, $|V_1\cap
    V_2|$ is equal to $[T^\wedge:M_1+M_2]$, the lattice index of $M_1+M_2$ in $T^\wedge$.
    \end{prop}

    \begin{proof}  The following argument is adapted from
    \cite[pp.32-33]{Systems}.
    Pick bases for $M_1$ and $M_2$.   $V_i$ is cut out by the equations
    $$x^{\mathbf a}=y_1^{\mathbf a},\ \ x^{\mathbf b}=y_2^{\mathbf b}$$
    for $x\in(\kay^*)^n$
    where ${\bf a}$ ranges over the basis for $M_1$ and ${\bf b}$ ranges over a basis for $M_2$.
    We write the basis vectors as row vectors and concatenate
    them to form an $n\times n$-matrix.
    $$A=\left[\begin{array}{c}
       A_1\\
       A_2
       \end{array}\right]$$
    Put this matrix in Hermitian normal form $UA=R$
    where $U\in \operatorname{SL}_n(\zee)$, and $R$ is an upper
    triangular invertible matrix.  Therefore, the entries of
    $R$ are
    $$R=\left[\begin{array}{cccc}
       r_{11}     &r_{12}   &\dots     &r_{1n}\\
       0          &r_{22}   &\dots     &r_{2n}\\
       \vdots     &\vdots   &          &\vdots\\
       0          &0        &\dots     &r_{nn}
       \end{array}\right].$$
    Finding intersection points of $V_1$ and $V_2$ amounts to
    solving the system
    \[x_1^{r_{i1}}x_2^{r_{i2}}\dots x_n^{r_{in}}=c_i\]
    for certain $c_i\in\kay^*$.
    There are $r_{11}r_{22}\dots
    r_{nn}=\det(A)=[T^\wedge:M_1+M_2]$ solutions.
    \end{proof}

    The definition of tropical intersection numbers in
    \cite{Mik06} requires that the above lattice index be
    equal to $[\zee^n:M_1^\perp+M_2^\perp]$ where $M_i^\perp$ is
    the perpendicular lattice to $M_i$.  For the sake of
    completeness, we include a proof with
    simplifications by Fr\'{e}d\'{e}ric Bihan that the lattice indexes are
    equal.

    \begin{lem} Let $L$ and $M$ be saturated lattices in
    $\zee^n$ of complementary rank so that $L+M$ has rank $n$.
    Then
    $$[\zee^n:L+M]=[\zee^n:L^\perp+M^\perp]$$
    where
    \[\begin{array}{rcl}
    L^\perp&=&\ker((\zee^n)^\vee\rightarrow L^\vee),\\
    M^\perp&=&\ker((\zee^n)^\vee\rightarrow M^\vee).
    \end{array}\]
    \end{lem}

    \begin{proof}
    Let $k=\rank(L)$.  Let $Q=\{q_1,\dots,q_k\}$ be a basis for
    $M^\perp$ and $R=\{r_1,\dots,r_k\}$ be a basis for $L$.

    We first claim that  
    \[[Z^n:L+M]=\left|\det\left([q_i(r_j)]_{i,j=1,\dots,k}\right)\right|.\] 
    Since $M$ is saturated, we may pick a basis $E=\{e_1,\dots,e_n\}$ for $\zee^n$ so that
    $\{e_{k+1},\dots,e_n\}$ is a basis for $M$. Let $F=\{f_1,\dots,f_k\}$ be a basis
    for $L$, and form the $n\times n$-matrix $A$ whose column
    vectors are the coordinates of $f_1,\dots,f_k,e_{k+1},\dots,e_n$
    with respect to the basis $E$.  $[\zee^n:L+M]=|\det(A)|$.
    The matrix $A$ is block lower-triangular with respect to
    blocks of size $k\times k$ and $(n-k)\times (n-k)$ centered at
    the diagonal.  The lower right $(n-k)\times (n-k)$ block is the
    identity matrix.  Therefore,
    $$\left|\det(A)\right|=\left|\det\left([a_{ij}]_{i,j=1,\dots k}\right)\right|=
    \left|\det\left([e_i^\vee(f_j)]_{i,j=1,\dots,k}\right)\right|.$$
    The determinant on the right is invariant under change of
    basis for $L$ and $M^\perp$.  The claim is proven.

    Similarly, $[\zee^n:L^\perp+M^\perp]$ is the absolute value of
    the determinant of the $k\times k$-matrix formed by letting a
    basis of $(L^\perp)^\perp$ act on a basis of $M^\perp$.  Since
    $L$ is saturated, $(L^\perp)^\perp=L$, so $R$ is a basis of
    $(L^\perp)^\perp$.  Therefore,
    $$[\zee^n:L^\perp+M^\perp]=\left|\det\left([r_i(q_j)]_{i,j=1,\dots,k}\right)\right|.$$
    It follows that the lattice indexes, $[\zee^n:L+M]$,$[\zee^n:L^\perp+M^\perp]$
    are equal to the absolute values of determinants of transposed matrices.  Therefore, they are equal.
    \end{proof}

    \subsection*{Definition of $\Trop$}
    Let $\cy$ be an immersive toric scheme defined over $\kay$ so $\cy=Y_0\times_{\Spec\ \kay} \Spec\ \ccr$.  Let $X$ be
    some subvariety of $\cy$ that intersects the big open torus.
    Let $\HI$ be the Hilbert image of $X$.  Its
    complex is the Gr\"{o}bner complex.

    \begin{defn}
    The tropical variety of $X$, $\Trop(X)\subset G^n$ is given by all
    $w\in G^n$ so that $\init{w}{X}$ intersects the big
    open torus, $(\kay^*)^n\subset Y_0$.
    \end{defn}

    By Proposition \ref{dualequivgrob}, if $w$ and $w'$ are in the same cell of the Gr\"{o}bner complex, then $\init{w}{X}$ is related to $\init{w'}{X}$ by an action of $(\kay^*)^n$.  So, if 
    $\init{w}{X}$ intersects the big open torus, so does $\init{w'}{X}$.   Therefore, the tropical variety
    is  a union of cells of the Gr\"{o}bner complex.   We may put a integral polyhedral complex structure on $\Trop(Y)$ to make it a subcomplex of the Gr\"{o}bner complex.
    
    The tropical variety is usually given by the image under the
    valuation map.  We show that these definitions are equivalent.

    Consider the isomorphism between the big open torus of $Y$ and
    $(\Kay^*)^n$ given by $g\mapsto g\cdot y$.
    This allows us to define a valuation map $v:X\cap (\Kay^*)^n\rightarrow G^n$

    \begin{lem} $\Trop(X)$ is equal to the image $-v(X)$.
    \end{lem}

    \begin{proof}
    $-v(X)\subseteq\Trop(X)$: Let $g\in X\cap(\Kay^*)^n$.  It suffices to show that the degeneration $\ol{g^{-1}\cdot X}\times_{\Spec\ \ccr}\Spec\ \kay$ intersects the big open torus in $Y_0$.  But,
    \[1=\ol{(g^{-1}\cdot g)}\times_{\Spec\ \ccr}\Spec\ \kay\in \ol{g^{-1}\cdot X}\times_{\Spec\ \ccr}\Spec\ \kay\]
    is a point in the big open torus.

    $\Trop(X)\subseteq -v(X)$: If $w\in\Trop(X)$, then
    \[\ol{t^w\cdot X}\times_{\Spec\ \ccr}\Spec\
    \kay\cap(\kay^*)^n\]
     is non-empty.  Let $\tilde{x}$ be a
    closed point of the above.  Then Lemma \ref{limitpoints}
    produces a point $x\in X$ with $\init{w}{x}=\tilde{x}$.  It follows that
    $-v(x)=w$.
    \end{proof}

   \begin{exmp} Let $H\subset T$ be a sub-torus and $x\in T$.  Let $X=H\cdot x$.  Then $\Trop(X)$ is $-H_G^\vee-v(x)$.
   \end{exmp}
   
   \begin{exmp} \label{trophypersurface} Let us revisit Example \ref{hypersurface}.  The Hilbert image is the toric variety associated to $-\ca$.  We have the morphism $(\proj^{|\ca|-1})^\vee\rightarrow\HI$.  The hypersurface in $Y$ corresponding to $[a_\omega]\in\HI$,
   \[\sum_{\omega\in\ca} a_\omega x^\omega=0\]
    is disjoint from the big open torus if and only exactly one $a_\omega$ is not zero.  Such points correspond to the torus fixed points of $\HI$ or alternatively, the top-dimensional cones of the Gr\"{o}bner fan.  Therefore the tropical variety of the hypersurface $V(f)$ is the union of the positive codimension cones of $N(-\Conv(\ca))$.
   \end{exmp}

    Let us relate the tropical variety of $\init{w}{X}$ to that of $X$.  
    
    \begin{lem} \label{tropicalstar} Let $w$ be a point in a cell $\tau$ of the tropical variety, $\Trop(X)$.  Then $\Trop(\init{w}{X})$ is the star of $\tau$ in $\Trop(X)$.
    \end{lem}
    
    \begin{proof}   
    Recall that by Lemma \ref{smalldef}, $\init{u}{{\init{w}{X}}}=\init{{w+\epsilon u}}{X}$ for sufficiently small $\epsilon$.  Therefore, 
    $\init{u}{{\init{w}{X}}}$ 
    intersects the open torus if and only $w+\epsilon u\in\Trop(X)$.   
     \end{proof}

    The dimension of $X$ and the dimension of $\Trop(X)$ are
    related.  We give a proof adapted from \cite{Systems}.
    We begin with the case where $\Trop(X)$ is
    zero-dimensional.

    \begin{lem}
    \label{zerodim}
    If $X\subseteq (\Kay^*)^n$ is a variety with
    $\dim(\Trop(X))=0$ then $X$ is zero-dimensional.
    \end{lem}

    \begin{proof}
    Suppose $X$ is positive dimensional.  Choose a coordinate projection
    $p:(\Kay^*)^n\rightarrow \Kay^*$ so that $p(X)$ is an
    infinite set.  By Chevalley's theorem \cite{RedBook}, $p(X)$
    is a finite union of locally closed sets and, since it is
    infinite, it must be an open set.  Therefore, $\Trop(X)$ is
    bigger than a point.
    \end{proof}

    We can reduce the general case to the above lemma.
    \begin{prop}
    \label{tropdim}
    If  $X\cap(\Kay^*)^n$ is purely
    $d$-dimensional, so is $\Trop(X)$.
    \end{prop}

    \begin{proof}
    Suppose $\dim\Trop(X)=k$.  Let $w$
    be an element of the relative interior of a top-dimensional
    cell of $\Trop(X)$.  Then $w$ is in the relative interior of a
    $k$-dimensional cell $C_\Gamma$ of the Gr\"{o}bner complex.  Then, by Lemma
    \ref{InvariantTorus}, $\init{w}{X}$ is invariant under a
    $k$-dimensional torus, $U$.   The initial degeneration $\init{w}{X}$ intersects the open torus so if $x\in\init{w}{X}\cap(\kay^*)^n$, the $k$-dimensional variety $U\cdot x$ is a subset of $\init{w}{X}$.  
    Since $\init{w}{X}$ is a flat
    deformation of $X$, it is also $d$-dimensional.  Therefore $k\leq d$.
    By Lemma \ref{tropicalstar}, the tropical variety of $\init{w}{X}$ is the $k$-dimensional subspace $\Span(C_\Gamma-w)$. 
    
    Now, we show $d=k$.  Let $W$ be a variety of the form $H\cdot z$ where $H\subset (\kay^*)^n$ is an $(n-k)$-dimensional torus with $H^\vee$ is transverse to 
    $\Trop(\init{w}{X})$.  Now, by the Kleiman-Bertini theorem
    \cite{KleimanBertini}, there is a choice of $z$ so that
    $\init{w}{X}\cap W$ is empty or of dimension $d-k$.  By Proposition \ref{intersectsubtori},
    $U\cdot x$ and $W$ must intersect, so $\init{w}{X}\cap W$ is non-empty.  
    But, $\Trop(\init{w}{X}\cap W)\subseteq \Trop(\init{w}{X})\cap\Trop(W)$ which is a point.
    Therefore, $\init{w}{X}\cap W$ is a $d-k$ dimensional scheme whose
    tropicalization is a point.  By the above lemma $d=k$.
    \end{proof}

    \subsection*{Multiplicities}

    Let $X$ be an
    $m$-dimensional subvariety of a toric variety $Y$.
    If $w$ is in the relative interior of an $m$-dimensional cell
    $C_\Gamma$ of $\Trop(X)$, then $\init{w}{X}\cap (\kay^*)^n$ is a
    subscheme invariant under an $m$-dimensional torus $H$ with $H_\real^\vee=\Span(C_\Gamma-w)$.
    Therefore, $\init{w}{X}\cap (\kay^*)^n$ is supported on
    $\coprod_i (H\cdot p_i)$ where $p_i$ are points in $(\kay^*)^n$.  This allows
    us to define multiplicities on $\Trop(X)$. 

    \begin{defn} Given a top-dimensional cell $\sigma$ of
    $\Trop(X)$, let $w$ be a point in the relative interior of
    $\sigma$.  Decompose the underlying cycle of $\init{w}{X}\cap (\kay^*)^n$ as
    $$[\init{w}{X}\cap (\kay^*)^n]=\sum m_i [H\cdot p_i]$$
    for $H\cong(\kay^*)^m\subset T_\kay$, $p_i\in (\kay^*)^n$.  The {\em multiplicity} $m_\sigma$ is
    $$m_\sigma=\sum_i m_i.$$
    \end{defn}

    This multiplicities are also called weights.

    $\Trop(X)$ obeys the following {\em balancing condition} first given in Theorem 2.5.1 of \cite{Speyer}.

    \begin{defn} An integrally weighted $m$-dimensional integral polyhedral complex is said
    to be {\em balanced} if the following holds: Let $\tau$ be an $(m-1)$-dimensional cell of $\Trop(X)$ and
    $\sigma_1,\dots,\sigma_l$ be the $m$-dimensional cells adjacent to $\tau$.  Let $w\in\tau^\circ$, $V=\Span(\tau-w)$, and $\lambda$ the projection $\lambda:T^\vee\rightarrow T^\vee/V$.
  Let $p_j=\lambda(\sigma_j-w)$.   Note that $p_j$ is an interval adjacent to $0$, and let $v_j\in T^\vee/V$ be the primitive integer vector along $\Span_+(p_j)$.  Then
    $$\sum_{j=1}^l m_{\sigma_j} v_j=0.$$
    \end{defn}

    We will give a proof that the balancing condition is satisfied in Theorem \ref{balancing}.  The following relates the multiplicities on $\Trop(\init{w}{X})$ to those on $\Trop(X)$.
    
    \begin{lem} \label{multd} Let $w\in\tau^\circ$ be a point in the relative interior of a cell of $\Trop(X)$.  Let $\sigma_1,\dots,\sigma_l$ be the top-dimensional cells in $\Trop(X)$ containing $\tau$.  Then the multiplicities of the cones $\ol{\sigma}_1,\dots,\ol{\sigma}_l$ in $\Trop(\init{w}{X})$ corresponding to $\sigma_1,\dots,\sigma_k$ are $m_{\sigma_1},\dots,m_{\sigma_l}$.
    \end{lem}
    
    \begin{proof}
      Let $u\in\ol{\sigma}_i$.  Then $\init{u}{\init{w}{X}}=\init{{w+\epsilon u}}{X}$ by Lemma \ref{smalld}.  By shrinking $\epsilon$ further if necessary, we may suppose $w+\epsilon u\in\sigma_i$.  Therefore, the degeneration $\init{u}{\init{w}{X}}$ used to compute $m_{\ol{\sigma}_i}$ is the same as the degeneration $\init{{w+\epsilon u}}{X}$ used to compute $m_{\sigma_i}$.
  \end{proof}

    \section{Intersection Theory Motivation: Bezout vs. Bernstein}

    Let us consider two curves in $(\com^*)^2$ cut out by
    polynomials $f(x,y)$ and $g(x,y)$.  Suppose they have no component in
    common and we would like to bound the number of intersection
    points in $(\com^*)^2$ counted with multiplicity.  The Bernstein bound will motivate tropical intersection theory.  

    \subsection*{Bezout Bound}
    We first consider the Bezout bound.  We compactify
    $(\com^*)^2$ to the projective plane $\proj^2$.  The
    intersection number is given by topology and is equal to
    $\deg(f)\deg(g)$.  This intersection bound is rigid in that
    it is invariant under deformations of $f$ and $g$.
    Unfortunately, the bound is not the best because we introduced
    new intersections on the coordinate hyperplanes by
    compactifying.

    Let us make this concrete by picking polynomials (all borrowed
    from \cite{Systems}).  Let
    \begin{eqnarray*}
    f(x,y)&=&a_1+a_2x+a_3xy+a_4y\\
    g(x,y)&=&b_1+b_2x^2y+b_3xy^2.
    \end{eqnarray*}
    To consider these polynomials on $\proj^2$, we must homogenize
    them to
    \begin{eqnarray*}
    F(X,Y,Z)&=&a_1Z^2+a_2XZ+a_3XY+a_4YZ\\
    G(X,Y,Z)&=&b_1Z^3+b_2X^2Y+b_3XY^2.
    \end{eqnarray*}
    Then the Bezout bound is $2\cdot 3=6$.  Notice that
    both curves contain the points $[1:0:0]$ and $[0:1:0]$.  This
    leads Bezout's theorem to over-count the number of intersections by $2$.
    It is impossible to remove these additional intersection points
    by an action of $(\com^*)^2$ since these points are fixed under the torus action.

    \subsection*{Bernstein Bound} Another approach is offered by Bernstein's theorem:
    \begin{thm}
    Given Laurent polynomials
    $$f_1,\dots,f_n\in\com[x_1^{\pm 1},\dots,x_n^{\pm 1}]$$
    with finitely many common zeroes in $(\com^*)^n$, let
    $\Delta_i$ be the Newton polytopes of $f_i$.  The number of
    common zeroes is bounded by the mixed volume of the
    $\Delta_i$'s.
    \end{thm}

    Bernstein's theorem can be conceptualized in the above case as follows.  One can
    compactify $(\com^*)^2$ to a nonsingular toric variety so that the closure of the curves cut
    out by $f=0$ and by $g=0$ does not intersect any torus fixed points.  For
    instance, one may take the toric variety whose fan is the
    normal fan to the Minkowski sum of the Newton polygons of $f$ and $g$.
    One may apply a $(\com^*)^2$-action to $\ol{\{f=0\}}$ to
    ensure that there are no intersections outside of
    $(\com^*)^2$.  By refining the fan further, we may suppose that the toric variety is smooth.
    Then one can bound the number of intersection points by the
    topological intersection number of the two curves.  This
    reproduces the Bernstein bound.

    \section{Intersection Theory}

    Henceforth, we will be using tropical varieties $Y(\Delta)$ defined by a fan $\Delta$ as in 
    \cite{FultonToric}.
    
    \subsection*{Intersection Theory over Discrete Valuation Rings}
    
    Let us first review some notions of Intersection Theory from \cite{Fulton}.  Let $Y$ be a scheme.  A $k$-cycle on $Y$ is a finite formal sum,  $\sum n_i [V_i]$
    where the $V_i$'s are $k$-dimensional subvarieties of $Y$ and the the $n_i$'s are integers.  $k$-cycles form a group under formal addition.  There is a notion of rational equivalence on cycles, and the Chow group, $A_k(Y)$ is the group of cycles defined up to rational equivalence.  This group is analogous to homology.    If $Y$ is complete, there is a natural degree map $\deg:A_0(Y)\rightarrow \zee$
    given by 
    \[\sum m_i[p_i]\mapsto \sum m_i.\]
    For any proper morphism $f:X\rightarrow Y$, there is an induced push-forward homomorphism
    \[f_*:A_k(X)\rightarrow A_k(Y).\] This push-forward homorphism commutes with degree.  If $X$ is a disjoint union $X=\bigsqcup X_i$, then we have 
    $A_k(X)=\bigoplus A_k(X_i)$.
       If $Y$ is a smooth $n$-dimensional variety, there is an intersection product 
\[A_k(Y)\otimes A_l(Y)\rightarrow A_{k+l-n}(Y).\]
If $V$ and $W$ are varieties in $Y$ of dimension $k$ and $l$ respectively, then the intersection product factors through a refined intersection product 
\[\xymatrix{A_k(Y)\otimes A_l(Y)\ar[r]& A_{k+l-n}(V\cap W)\ar[r]^>>>>>>{i_*}&A_{k+l-n}(Y)}\]
where $i:V\cap W\rightarrow Y$.   
There is also Chow cohomology $A^k(Y)$ which is defined operationally.

    Intersection theory can also be defined over discrete valuation rings.  The reference is Chapter 20 of \cite{Fulton}.  We will state the results for
    $\ccr=\com[[t^{\frac{1}{M}}]]$, but they are true for more general
    choices of $\ccr$.  In practice, however, given varieties defined over $\com\{\{t\}\}$, we may find a
    sufficiently large $M$ so that they are defined over $\com((t^{\frac{1}{M}}))$ and apply the results
    for the corresponding choice of $\ccr$.
    Let $p:\cy\rightarrow\Spec\ \ccr$ be a scheme over $\Spec\ \ccr$.
    Let $Y=\cy\times_{\Spec\ \ccr}\Spec\ \Kay$, $Y_0=\cy\times_{\Spec\ \ccr}\Spec\
    \kay$.

    Many results from intersection theory including the existence of degree and refined intersection product remain true in this case using relative dimension
    over $\Spec\ \ccr$ in place of absolute dimension.  The new feature in this situation is the {\em specialization
    map}
    \[s:A_k(Y/\Kay)\rightarrow A_k(Y_0/\kay)\]
    which is the Chow-theoretic analog of
    $X\rightarrow (\ol{X})\times_{\Spec\ \ccr}\Spec\ \kay$.

    \begin{prop} If $\cy$ is smooth over $\Spec\ \ccr$ then
    the specialization map is a ring homomorphism.  Moreover it
    commutes with refined intersection product.
    \end{prop}

    \begin{proof}
    See Corollary 20.3  and Example 20.3.2 in \cite{Fulton}.
    \end{proof}
    \subsection*{Transversal Intersections}

    Let $V^k,W^l\subset Y^n$ be varieties of dimensions $k$
    and $l$ where $k+l=n$.  Let $Y$ be a smooth toric variety over
    $\Spec\ \Kay$.

    \begin{defn} $V^k$ and $W^l$ are said to {\em intersect properly} if $V\times_Y W$ is a
    zero-dimensional scheme.
    \end{defn}

    \begin{defn} Two tropical varieties $\Trop(V)$,
    $\Trop(W)$ are said to {\em intersect transversally} if
    they intersect in the relative interior of transversal top-dimensional cells.
    \end{defn}

    Note that it is not sufficient that $V$ and $W$ intersect
    transversally for $\Trop(V)$ and $\Trop(W)$ to intersect
    transversally.  In fact, $V$ and $W$ can be disjoint
    while their tropicalizations intersect (or even coincide, for example, $x+y=1$ and $x+y=2$ in $(\Kay^*)^2$.
    However, the transversal intersection lemma of \cite{Computing} does give a
    condition for $V$ and $W$ to intersect:

    \begin{lem} If $\Trop(V)$ and $\Trop(W)$ intersect
    transversally at $w\in\real^n$, then $w\in\Trop(V\cap W)$.
    \end{lem}

    \begin{proof}
    Since $w$ is in a top dimensional cell of $\Trop(V)$ and of
    $\Trop(W)$ then
    \[\supp(\init{w}{V})=H_1\cdot V_\sigma\]
    \[\supp(\init{w}{W})=H_2\cdot W_\tau\]
    where $\supp$ denotes underlying sets, $V_\sigma,W_\tau$ are finite sets of points, and
    $H_1,H_2$ are sub-tori of dimension $k$ and $l$, respectively.    By Proposition
    \ref{intersectsubtori}, 
    \[(\init{w}{V}\times_{Y_0}\init{w}{W})\cap (\kay^*)^n\]
     is non-empty and
    zero-dimensional.  Let $z$ be a closed point of $(\init{w}{V}\times_{Y_0}\init{w}{W})\cap
    (\kay^*)^n$.  Now let $\ccv=\ol{t^w\cdot V}, \cw=\ol{t^w\cdot W}$.  Let $\cz$
    be a maximal irreducible component of
    $\ccv\times_{\cy} \cw$ containing $z$.  Therefore,
    $(\cz\times_{\Spec\ \ccr}\Spec\ \kay)\cap(\kay^*)^n$ is non-empty and zero-dimensional.  We
    claim $\cz$ is not contained in the fiber over $\Spec\ \kay$.

    \begin{claim} $\cz$ surjects onto $\Spec\ \ccr$.
    \end{claim}

    Since $\ccv$ and $\cw$ have relative dimension $k$ and $l$,
    respectively, each top-dimensional irreducible component $\ccv\times_{\cy}\cw$ must have relative
    dimension at least $0$ and therefore cannot be contained in the special fiber
    as a $0$-dimensional subscheme.

    $Z=\cz\times_{\Spec\
    \ccr}\Spec\ \Kay\subset t^w V\times_Y t^w W$ is non-empty and
    $z\in\init{w}{t^{-w}Z}\subseteq \init{w}{V\times_Y W}$.  Therefore $V\times_Y W$ must
    have a point of valuation $-w$.
    \end{proof}

    \begin{lem} If all intersections of $\Trop(V)$ and
    $\Trop(W)$ are transversal, then $V\cap (\Kay^*)^n$ and $W\cap (\Kay^*)^n$
    intersect properly.
    \end{lem}

    \begin{proof}
    Let $Z$ be the intersection of the two varieties with the
    reduced induced structure.  Then
    $\Trop(Z)=\Trop(V)\cap\Trop(W)$ is zero-dimensional.
    Lemma \ref{zerodim} shows that every component of $Z$ is
    zero dimensional.
    \end{proof}

    \subsection*{Intersection of Tropicalizations}

   We will define an intersection number for transversal tropical varieties of complementary dimensions.

    Let $Y$ be an $n$-dimensional smooth toric variety defined over $\kay$.
    Let $V^k,W^l\subseteq Y$ be varieties of
    complementary dimensions such that $\Trop(V)$ and $\Trop(W)$
    intersect tropically transversely.  Let
    $x\in\Trop(V)\cap\Trop(W)$ such that $x$ is contained in top-dimensional cells $\sigma_x,\tau_x$ of $\Trop(V)$ and $\Trop(W)$, respectively.    Translate $\Trop(V)$ and $\Trop(W)$ so that $x$ is at the origin.  We have inclusions
    $\real\sigma_x, \real\tau_x\hookrightarrow T^\vee_\real$ which induce
    projections
    $T^\wedge_\real\rightarrow (\real\sigma_x)^\vee$ and
    $T^\wedge_\real\rightarrow (\real\tau_x)^\vee$.
     Let $M_x$ and $N_x$ be the lattices defined by
    \[M_x=\ker(T^\wedge_\real\rightarrow (\real\sigma_x)^\vee)\cap T^\wedge\]
    \[N_x=\ker(T^\wedge_\real\rightarrow (\real\tau_x)^\vee)\cap T^\wedge\]
    Let $m_x, n_x$ be the multiplicities of $\sigma_x$ and
    $\tau_x$ in $\Trop(V)$ and $\Trop(W)$ respectively and define
    the {\em tropical intersection number} to be
    \[\deg(\Trop(V)\cdot\Trop(W))=\sum_{x\in\Trop(V)\cap\Trop(W)} m_x n_x[T^\wedge:M_x+N_x].\]
    This definition is analogous to the definition in classical intersection theory.  Here, $m_x,n_x$ are analogous to the multiplicities of subvarieties in cycles and the lattice index is analogous to the length of a zero-dimensional component of the intersection.  

   \begin{defn} $V$ and $W$ {\em intersect in the interior} if
    the support of $V\times_Y W$ is contained in the big open torus $T$ of $Y$.
    \end{defn}

    \begin{thm}
    \label{IntersectionNumber}
    If $V$ and $W$ intersect tropically
    transversally and in the interior then the tropical intersection
    number of $\Trop(V)$ and $\Trop(W)$is equal to the classical intersection number.
    \end{thm}

    \begin{proof}
    Let us replace $\Kay$ by a field $\com((t^{\frac{1}{M}}))$ over which $V$
    and $W$ are defined.
    First note that $\Trop(V\cap
    W)=\Trop(V)\cap \Trop(W)$ by the transverse intersection lemma.  It follows that $V\cap W$ is zero-dimensional.
    Decompose this intersection into a disjoint union
    $$V\times_Y W = \coprod_{x\in\Trop(V)\cap\Trop(W)} Z_x$$
    where $v(Z_x)=-x$.
    Now, the refined intersection product is
    \[V\cdot W\in A_0(V\cap W)=\bigoplus A_0(Z_x)\]
    and the intersection number is the degree of the intersection
    product.  If 
    \[\pi_x:A_0(V\cap W)\rightarrow A_0(Z_x)\]
     is the projection onto the summand, then 
    \[\deg(V\cdot W)=\sum_{x\in\Trop(V)\cap\Trop(W)} \deg(\pi_x(V\cdot W)).\]
    
    Let $w\in \Trop(V)\cap\Trop(W)$ 
    and
    \[\ccv=\ol{t^w\cdot V}\subseteq \cy\]
    \[\cw=\ol{t^w\cdot W}\subseteq \cy.\]
    Note that $\ccv$ and $\cw$ are flat over $\ccr$.

    Decompose the intersection of $\ccv$ and $\cw$ as
    \[\ccv\times_{\cy} \cw = \coprod_{x\in\Trop(V)\cap\Trop(W)} \cz_x\]
    where
    \[\cz_x\times_{\Spec\ \ccr}{\Spec\ \Kay}=t^w\cdot Z_x.\]
    The zero-dimensional scheme $(\cz_x)_0=\cz_x\times_{\Spec\ \ccr}\Spec\ \kay$ is contained in $(\kay^*)^n$ only if $x=w$.  Otherwise, it is disjoint from $(\kay^*)^n$.   
    Let   $(\ccv\times_{\cy}\cw)_0=(\ccv\times_{\cy}\cw)\times_{\Spec\ \ccr}\Spec\ \kay$.
    Since $\cz_w$ is proper over $\Spec\ \ccr$, by Proposition 20.3 and Corollary 20.3 of \cite{Fulton}, the image of $[t^w V]\otimes [t^w Y]$ under the following compositions are equal
    \[\xymatrix{A_k(Y)\otimes A_l(Y)\ar[r]&A_0(t^w(V\cap W))\ar[r]^>>>>>>>>>>{\pi_w}&A_0(t^wZ_w)\ar[r]^{s}&A_0((\cz_w)_0)\ar[r]^>>>>>>>{\deg}&\zee\\
                      A_k(Y)\otimes A_l(Y)\ar[r]^{s\otimes s}&A_k(Y_0)\otimes A_l(Y_0)\ar[r]&
                      A_0((\ccv\times_\cy\cw)_0)\ar[r]^>>>>>>{\pi_w}&
                      A_0((\cz_w)_0)\ar[r]^>>>>>>>{\deg}&\zee.}\]
    But the second composition is just the degree of the intersection of 
    the tori $\init{w}{V}$ and
    $\init{w}{W}$.  Their intersection number is $m_wn_w[T^\wedge:M_w+N_w]$ by Proposition \ref{intersectsubtori}.  Summing over $w\in\Trop(V)\cap\Trop(W)$, we get the
    result.
    \end{proof}

    \subsection*{Transversality}

    \begin{lem}
    \label{ConstantTransversality}
    If $V$ and $W$ intersect all torus orbits properly then there exists $\lambda\in(\kay^*)^n$, such that
    $\lambda\cdot V$ intersects $W$ properly and in the interior.
    \end{lem}

    \begin{proof}
    By the Kleiman-Bertini theorem \cite{KleimanBertini} applied to each orbit closure $V(\sigma)$, 
    there exists a non-empty open set $U\subset (\Kay^*)^n$ such that for all $\lambda\in
    U$, $\lambda\cdot V$ intersects $W$ properly and in the interior.  It
    suffices to show that $U\cap (\kay^*)^n$ is non-empty.

    Suppose $U\cap (\kay^*)^n$ is empty.  Let $f\in\Kay[x_1^{\pm 1},\dots,x_n^{\pm 1}]$ be a Laurent
    polynomial over $\Kay$ so that $(\Kay^*)^n\setminus V(f)\subseteq U$.  Then $V(f)$
    contains all $\kay$-points.  By clearing denominators, we may
    suppose $f\in R=\com[[t^{\frac{1}{M}}]]$ for some $M$ where $t^{\frac{1}{M}}$
    does not divide $f$.  Since $f=0$ on $(\kay^*)^n$, $f|_{t^{1/M}=0}=0$.  It follows that $t^{\frac{1}{M}}$ divides $f$.  This gives a contradiction.
    \end{proof}

    Note that $\lambda\cdot V$ and $V$ have the same tropical
    variety.

    \subsection*{Balancing Condition}

    In this section, we prove that if $X$ is an $m$-dimensional subvariety of a
    toric variety $Y$, then $\Trop(X)$ satisfies the balancing
    condition.  The strategy of the proof is that a well-defined tropical intersection number between $\Trop(X)$ and $\Trop(H\cdot z)$ for $H$ a sub-torus and $z\in T$ guarantees that $\Trop(X)$ is balanced.

     We need the following technical lemma.
     
     \begin{lem} Let $x\in (\Kay^*)^n$ and $\sigma$ be a cone in $\Delta$. Then $v(x)\in
\sigma^\circ$ if and only if $\init{0}{x}\in\oh_\sigma$, the open torus corresponding to $\sigma$.  
\label{toruslimits}
\end{lem}

\begin{proof}
Consider the toric chart $U_\sigma=\Spec\ \Kay[\sigma^\vee\cap T^\wedge]\supset(\Kay^*)^n$. The torus orbit
$\oh_\sigma$ is cut out by the ideal $I_\sigma$ which is the kernel
of the projection
\[\Kay[\sigma^\vee\cap T^\wedge]\rightarrow\Kay[\sigma^\perp\cap T^\wedge].\]
A monomial $m\in I_\sigma$ is of the form $x^u$ for $u$ satisfying $\<u,y\>>0$ for all $y\in\sigma^\circ$.
Since $v(x)\in \sigma^\circ$, $v(m(x))>0$ for every monomial $m\in
I_\sigma$ while $v(m(x))=0$ for every $m\in\Kay[\sigma^\perp\cap T^\wedge]$.  

Suppose $v(x)\in\sigma^\circ$.
If $m\in\com[\sigma^\vee\cap T^\wedge]$ is a monomial, $m(x)|_{t=0}=m(\init{0}{x})$.  Therefore, for $f\in I_\sigma$,
$v(f(x))>0$ so under the specialization $t=0$, $f(x)$ goes to $0$.  On the other hand, for every $m\in\com[\sigma^\perp\cap T^\wedge]$, $m(x)$ goes to its leading term which is non-zero.  It follows that $\init{0}{x}\in \oh_\sigma$.  

Now, suppose $\init{0}{x}\in\oh_\sigma$.  For any monomial $m=x^u\in I_\sigma$, we have $m(x)|_{t=0}=m(\init{0}{x})=0$.  Therefore, $v(m(x))>0$ which implies $\<u,v(x)\>>0$.  For $u\in\sigma^\perp$, $m=x^u$ is non-zero on $\init{0}{x}$.  It follows that $\<u,v(x)\>=v(m(x))=0$ and so $v(x)\in\sigma^\circ.$ 
\end{proof}

    We need the following Lemma of Tevelev.
    
    \begin{lem} \cite[Lemma~2.2]{Tevelev}
    \label{FromTevelev}
    Let $Y(\Delta)$ be a complete toric variety given by a fan $\Delta$.  Let
    $X\subset Y(\Delta)$ be a subvariety defined over $\kay$.  
    Then $-\Trop(X)$ intersects a cone $\sigma$ in the fan $\Delta$
    in its relative interior if and only if $\ol{X}$ intersects $\oh_\sigma$.
    \end{lem}

   \begin{proof}
Write $X_\kay$ for $X$ and $X_\Kay$ for $X\times_{\Spec\ \kay}\Spec\ \Kay$.
Observe that $X_\kay=\init{0}{X_\Kay}$.

Suppose $-\Trop(X)\cap\sigma^\circ$ is non-empty.  Then there
exists $x\in X_\Kay$ with $v(x)\in \sigma^\circ$.  Therefore,
$\init{0}{x}\in\oh_\sigma$.

Now suppose $\ol{X}\cap \oh_\sigma$ is non-empty.  Then by Corollary
\ref{openlimitpoints}, there exists $x\in V_\Kay\cap (\Kay^*)^n$ with
$\init{0}{x}\in \oh_\sigma$.  It follows that
$v(x)\in\sigma^\circ$.
\end{proof}

    \begin{defn} A subvariety $X\subset Y$ of dimension $l$ is said
    to {\em intersects orbits properly} if
    \begin{enumerate}
    \item[(1)] For $\sigma$ a cone in $\Delta$ with $\dim\sigma>l$,  $X\cap
    \oh_\sigma=0$.
    \item[(2)] For $\sigma$ a cone in $\Delta$ with
    $\dim\sigma=l$, $X\cap\oh_\sigma$ is a $0$-dimensional scheme.
    \end{enumerate}
    \end{defn}
  
By replacing $\Delta$ with a finer fan so that $-\Trop(X)$ is supported on a union of cones of dimension at most $l$, we may always ensure that $X$ intersects orbits properly.

We first prove that curves defined over $\kay$ are balanced.

\begin{lem} Let $X$ be a curve defined over $\kay$ in a complete toric variety $Y(\Delta)$.  Then $\Trop(X)$ is balanced.
\end{lem}  

\begin{proof}
By refining $\Delta$, we may suppose that $X$ intersects torus orbits properly and that $Y$ is smooth.
$\Trop(X)$ consists of rays $\rho_1,\dots,\rho_l$ weighted with multiplicities $m_1,\dots,m_l$.  Let $v_i$ be the primitive integer vector along $\rho_i$.
It suffices to show that 
\[\sum_{j=1}^l m_j \<u,v_j\>=0\]
for any $u$ in $T^\wedge$.
Let $H\subset T$ be the sub-torus so that $H^\vee=u^\perp$.  Let $W_y=H\cdot y$ for $y\in(\kay^*)^n$.
By refining $\Delta$ further, we may suppose that $W_y$ intersects torus orbits properly.  By replacing $W_y$ by $\lambda\cdot W_y$, we may suppose that $W_y$ intersects $X$ is the interior.  

Since for $w,w'\in T_G^\vee$, $t^wW$ and $t^{w'}W$ are related by the $T$-action, they are linearly equivalent. Therefore, by Lemma
    \ref{ConstantTransversality} and Theorem \ref{IntersectionNumber}, the tropical intersection number $\deg(\Trop(X)\cdot \Trop(t^wW))$ is
    independent of $w$.
    
    We may suppose without loss of generality that $u$ is primitive.  Pick $w\in T^\vee$ such that $\<u,w\> > 0$ and $y\in (\kay^*)^n$.  Then $\Trop(t^wW_y)=-w-H_\real^\vee$ with some multiplicity $n_W$.
Then $\rho_j\cap\Trop(t^w W)$ is non-empty if and only if $\<u,v_j\><0$.  The multiplicity of such an intersection is
\[m_jn_W[T^\wedge:(\zee u)+v_j^\perp]=m_jn_W|\<u,v_j\>|.\]
Therefore, 
\[\deg(\Trop(W)\cdot \Trop(X))=\sum_{j:\<u,v_j\><0} m_jn_W|\<u,v_j\>|.\]
Replacing $w$ by $-w$, we see
\[-\sum_{j:\<u,v_j\><0} m_jn_W\<u,v_j\>=\sum_{j:\<u,v_j\>>0} m_jn_W\<u,v_j\>\]
from which the conclusion follows.
\end{proof}

    \begin{thm} \label{balancing} $\Trop(X)$ satisfies the balancing condition.
    \end{thm}

    \begin{proof}
    Let $\tau$ be some $(m-1)$-dimensional cell of $\Trop(X)$ and
    $\sigma_1,\dots,\sigma_l$, the adjacent $m$-dimensional cells.  Let $w$ be a
    point in the relative interior of $\tau$. 
    $\init{w}{X}$ is a subscheme that
    is invariant under an $(m-1)$-dimensional torus.  $\Trop(\init{w}{X})$, the star of $\tau$
    consists of the linear subspace $\overline{\tau}=\Span(\tau-w)$ and the cones $\overline{\sigma_i}=\Span^+(\sigma_i-w)+\overline{\tau}$.  The multiplicities
    of the $\sigma$'s in $\Trop(\init{w}{X})$ are that same as those
    of the corresponding cells in $\Trop(X)$ by Lemma \ref{multd}.  
    
     Let $V$ be the
    union of the components of $\init{w}{X}$ that intersect the big
    open torus.  Then, $\Trop(V)=\Trop(\init{w}{X})$ and by refining $\Delta$, we may ensure 
    $V$ intersects the torus orbits properly.  Let $K$ be the
    $(m-1)$-dimensional invariant torus of $V$, and $p:T\rightarrow
    T/K$ be the quotient map.  
    The image of $\Trop(V)$ under that map is a one-dimensional
    integral polyhedral complex with one vertex and $l$ rays
    $\real_+ v_1',\dots,\real_+ v_l'$ emanating from it where $v_i'$ is a primitive integer vector.
    For $u\in(T/K)^\wedge$, let $H\subset (T/K)^\vee$ be the $(n-m-1)$-dimensional torus with $H^\vee=u^\perp$.  Now let $H'\subset T$ be a $(n-m-1)$-dimensional torus with $p(H')=H$.  Pick $w\in T^\vee_G$ such that $\<u,p^\vee(w)\>>0$.  For $y\in (\kay^*)^n$, let $W_y=H'\cdot y$. 
    Then $\ol{\sigma}_j$ intersects $\Trop(t^w W_y)$ if and only if $\<u,v_j\><0$.  The intersection multiplicity in that case is 
    \[\begin{array}{rcl}
    m_{\ol{\sigma}_j}n_W[T^\wedge:(H')^\wedge+(\ker(T^\wedge\rightarrow (\real\ol{\sigma}_j)^\vee)\cap T^\wedge)]&=&m_{\ol{\sigma}_j} n_W[(T/K)^\wedge:\zee u+v_j^\perp]\\
    &=&m_{\ol{\sigma}_j}n_W|\<u,v_j\>| .
    \end{array}\]
   The argument now proceeds as in the case of curves. 
    \end{proof}

   We should mention that the above argument can be simplified by using the theorem that tropical varieties are natural under monomial morphisms as proved by Sturmfels and Tevelev \cite{ST07}.

    \section{Tropical Cycles and the Cohomology of Toric Varieties}

    In this section, we work over a field $\Kay\supset\kay=\com$.  $\Kay$ may be the
    field of the Puiseux series or the complex numbers.

    \subsection*{Minkowski Weights}

    In \cite{FS}, Fulton and Sturmfels gave a description of Chow
    cohomology of a complete toric variety in terms of the fan.  This
    description is closely related to the balancing condition for
    tropical varieties.

    Consider a complete toric variety $Y$ given by a complete $n$-dimensional fan
    $\Delta$.  The Chow cohomology of
    $Y$ is given by Minkowski weights.  Let $\Delta^{(k)}$ be the
    set of all cones of codimension $k$.  For a cone
    $\sigma\in\Delta^{(k)}$, $\tau\in\Delta^{(k+1)}$, $\tau\subset
    \sigma$, let $N_\sigma$ be the lattice span of $\sigma$ and let
    $n_{\sigma,\tau}\in\sigma$ be an integer vector whose image generates the one-dimensional
    lattice $N_\sigma/N_\tau$.

    \begin{defn} A {\em Minkowski weight of
    codimension $k$} is a function
    $$c:\Delta^{(k)}\rightarrow \zee$$
    so that for every $\tau\in\Delta^{(k+1)}$ and every element $u\in \tau^\perp\cap
    \zee^n$,
    $$\sum_{\sigma\in\Delta^{(k)}|\sigma\supset\tau}
    c(\sigma)<u,n_{\sigma,\tau}>=0.$$
    \end{defn}

    As a consequence of showing $A^k(Y)=\Hom(A_k(Y),\zee)$, it is proven in  \cite{FS} that
    the Chow cohomology group $A^k(Y)$ is
    canonically isomorphic to the group of Minkowski weights of
    codimension $k$.

    We can view a Minkowski weight as an integrally  weighted integral fan,
    \[\bigcup_{c(\sigma)\neq 0} \sigma\]
    where the cone $\sigma$ is weighted by $c(\sigma)$.
    There is a formula for the cup-product in terms of Minkowski
    weights.  If we view Minkowski weights $c$ and $d$ of complementary dimension as fans, then their tropical intersection number (after translating one fan to ensure that they are tropically transverse) is equal to the degree of their cup product evaluated on the fundamental class of $Y,$ $\deg((c\cup d)\cap [Y])$.  
    
    If $X\subset Y$ is a codimension $k$ subvariety defined over $\kay$, the
    function taking a cone in $\Trop(X)$ to its multiplicity
    satisfies the balancing condition which is exactly the
    Minkowski weight condition.

    \subsection*{Associated Cocycles}
    If $Y$ is smooth, to every algebraic cycle $X$ of codimension $k$ in $Y$, we may associate
    a Minkowski weight of codimension $k$ by Poincare duality.  We will do this explicitly using toric geometry.
    
    \begin{lem} \label{poindual} Let $Y(\Delta)$ be a smooth toric variety over $\kay$.  Let $X$ be a codimension $k$ algebraic cycle.  Define a function 
\[\begin{array}{ccccl}  
 c&:&\Delta^{(k)}&\rightarrow&\zee\\ 
 c&:&\sigma&\mapsto& \deg([X]\cdot [V(\sigma)]).
  \end{array}\]        
Then $c$ is a Minkowski weight and $c\cap [Y(\Delta)]=[X]$.
\end{lem}

\begin{proof} $[X]$ has a Poincar\'{e}-dual $d$ satisfying $d\cap \alpha=\deg([X]\cdot \alpha)$ for $\alpha\in A_k(Y)$.
For all $k$-dimensional torus orbits, $V(\sigma)$, we have 
\[c(\sigma)=\deg([X]\cdot [V(\sigma)])=d\cap V(\sigma)\]
 Since $A_*(Y)$ is generated by torus orbits and $A^*(Y)=\Hom(A_*(Y),\zee)$, $c=d$ as Minkowski weights.
\end{proof}
    
    If $X$ is a subvariety of $Y$ defined over $\kay$, we may relax the smoothness condition on $Y$ after mandating that $X$ intersects the torus orbits of $Y$ properly.
    
    \begin{defn} Let $Y$ be a complete toric variety.  Let $\widetilde{Y}$ be a smooth toric resolution of $Y$  with fan $\widetilde{\Delta}$, which is a refinement of $\Delta$.
    Define the {\em associated cocycle of $X$}, a Minkowski weight on $\widetilde{\Delta}$ by $c(\tilde{\tau})=\deg([X]\cdot [V(\widetilde{\tau})])$.  
    \end{defn}
   
   The associated cocycle is well-defined as a Minkowski weight on $\widetilde{\Delta}$.  The following proposition shows that it is well-defined on $\Delta$.

    \begin{prop} If $X$ is an $k$-dimensional subvariety of $Y$, defined
    over $\kay$ that intersects the torus orbits properly then
    the associated cocycle of $X$ is $-\Trop(X)$.
    \end{prop}

    \begin{proof}
    Because $X$ intersects the torus orbits properly, by Lemma \ref{FromTevelev}, $-\Trop(X)$ is supported on
    $k$-dimensional cones in $\Delta$.

    We need only
    show that for every $\tilde{\tau}\in \Delta^{(k)}$, the multiplicity
    $m_{\tilde{\tau}}$ is equal to $c(\tilde{\tau})$.  Let $w\in -\tilde{\tau}^\circ$.
    Because intersection product commutes with specialization,
    $$\deg([X]\cdot [V(\tilde{\tau})])=\deg([\init{w}{X}]\cdot [V(\tilde{\tau})]).$$
    Let $H\subset T$ be the $k$-dimensional sub-torus corresponding to $\tau\subset T^\vee_\real$.
    The underlying cycle of $\init{w}{X}$ can be decomposed as
    $$[\init{w}{X}]=\sum m_i[H\cdot p_i]+D$$
    where $p_i\in(\kay^*)^n$ and $D$ is disjoint from the big open torus.

    We claim that $D$ is disjoint from $V(\tilde{\tau})$.  If it was not,
    it would have to intersect  a proper torus orbit of $V(\tilde{\tau})$.
    Therefore, it suffices to show that $\init{w}{X}$ does not intersect
    $V(\tilde{\sigma})$ for $\tilde{\sigma}\supset \tilde{\tau}$.  If it did, then by Corollary \ref{openlimitpoints},
    there would
    be $x\in X\cap (\Kay^*)^n$ so that $\init{w}{x}\in V(\tilde{\sigma})$.  By Lemma \ref{toruslimits},
    $v(x)+w\in \tilde{\sigma}^\circ$.  Therefore, $v(x)\in -w+\tilde{\sigma}^\circ\subset \tilde{\tau}^\circ+\tilde{\sigma}^\circ\subset\tilde{\sigma}^\circ$.
    But we assumed that $-\Trop(X)$ does not intersect $\tilde{\sigma}^\circ$ which is a cone
    of $\widetilde{\Delta}$ of dimension greater than $k$.

    By a local computation, we see $H\cdot p_i$ meets $V(\tilde{\tau})$ transversely in a single point.
    Therefore, $c(\tilde{\tau})=\sum m_i [H\cdot p_i]\cdot [V(\tilde{\tau})]=\sum m_i=m_{\tilde{\tau}}$.
    \end{proof}

   It follows that the associated cocycle is a pullback by $\pi:Y(\widetilde{\Delta})\rightarrow Y(\Delta)$.  Furthermore, the associated cocycle is dual to $[X]$.
   
   \begin{lem} If $c$ is the associated cocycle of $X\subset Y$, then
   \[c\cap [Y]=[X]\in A_k(Y)\].
    \end{lem}
    
    \begin{proof}
    Let $\pi:Y(\widetilde{\Delta})\rightarrow Y(\Delta)$ be a smooth toric resolution.  By Lemma \ref{poindual}, $\pi^*c\cap [Y(\widetilde{\Delta})]=[\pi^{-1}(X)]$. The projection formula tells us 
    \[c\cap [Y]=c\cap \pi_*([Y(\widetilde{\Delta})]=\pi_*([\pi^{-1}(X)])=[X].\]
    \end{proof}

    \begin{exmp} \label{multhypersurface}
    This gives us the weights for the tropicalization of the hypersurface found in Example \ref{trophypersurface}.  A top-dimensional cone of $\Trop(V(f))$ corresponds to a $1$-dimensional face $\Gamma\subset\Conv(-\ca)$.   The multiplicity of that cell is  $\deg(V(f)\cdot Y(\Gamma))$.  This intersection is defined by 
    \[\sum_{\omega\in\Gamma} a_\omega x^\omega=0.\]
    This is a polynomial in one variable whose Newton polytope is $\Gamma$.  Therefore, the number of points in the intersection, hence the multiplicity is the lattice length of the edge $\Gamma$.  
     \end{exmp}
    
    \begin{thm} Given two varieties $V^k$,
    $W^l$ that intersect torus orbits properly and intersect tropically transversely, the intersection number of their associated cocycles is equal to their tropical
    intersection number.
    \end{thm}

    \begin{proof}
    We pass to a smooth toric resolution.  By using
    the Kleiman-Bertini theorem, we may find $\lambda\in T_\kay$ so that $\lambda\cdot V$ and
    $W$ intersect in the interior. Note that
    $\Trop(\lambda\cdot V)=\Trop(V)$.  The intersection number of the
    associated tropical cycles is equal to the intersection pairing on their
    Poincare-duals in cohomology by \cite{FS}.  But this is their classical
    intersection number which equals
    $\deg(\Trop(V)\cdot\Trop(W))$ by Theorem \ref{IntersectionNumber}
    \end{proof}

    \subsection*{Proof of Bernstein's Theorem}

    For the sake of completeness, we outline a proof of Bernstein's theorem along
    the lines of the above section.  In essence, this proof is a hybrid of the
    proofs given in \cite{FultonToric} and \cite{Systems}.  We work over $\com$.

    Given Laurent polynomials
    \[f_1,\dots,f_n\in\com[x_1^{\pm 1},\dots,x_n^{\pm 1}],\]
    let $Q_i$ be the Newton polytope of $f_i$.  We summarize the facts we have established in the lemma below.

    \begin{lem} Let $f_i$ be a polynomial with Newton
    polytope $Q_i$, and $X(\Delta_i)$, the toric variety whose fan is $\Delta_i=N(Q_i)$  
    The hypersurface $V(f_i)$ intersects torus orbits in $X(\Delta_i)$ properly.
    \end{lem}

   We know by Example \ref{multhypersurface} that the tropical cycle $c_i$ associated to
    $V(f_i)$ is the union of cones of the normal fan of $\Delta_i$ of
    positive codimension where the codimension $1$ cones are weighted by the lattice
    length of the dual edges of $\Delta_i$.

    Let $\Delta$ be a fan that refines the
    normal fans of the $\Delta_i$'s so that $X(\Delta)$ is smooth.  There are birational morphisms
    from a nonsingular variety, $p_i:X(\Delta)\rightarrow X(\Delta_i)$.
    By \cite{Systems}, the mixed volume of
    $\Delta_1,\dots,\Delta_n$ is equal to the tropical intersection of the cycles
    $c_i$.  By \cite{FS}, this is equal to $\deg(p_1^*c_1\cup\dots\cup
    p_n^*c_n)$, which is the intersection number of
    $p_1^{-1}(V(f_1))\cdot\ldots\cdot p_n^{-1}(V(f_n))$ in $X(\Delta)$.  This bounds the number of geometric
    intersections in $(\com^*)^n$.

    \section{Deformations of Subschemes into Torus Orbits}

    This section is a generalization of Theorem 2.2  of \cite{DFS}.
    Let $Y(\Delta)$ be a smooth toric scheme defined over $\kay$ and $X\subseteq Y$, a purely $k$-dimensional
    closed subscheme. If $w$ is in the relative interior of an
    $m$-dimensional cell of the Gr\"{o}bner complex of $X$, then
    $\init{w}{X}$ is invariant under an $m$-dimensional torus.
    $\init{w}{X}$ has components supported in the big open torus of $Y$
    and within smaller dimensional torus orbits.  In particular if
    $w$ is in the interior of an open cell of the Gr\"{o}bner complex, $\init{w}{X}$ is invariant under $T$.
    Therefore, the maximal components of $\init{w}{X}$ are
    supported on the $k$-dimensional torus orbits.  We can use
    tropical geometry to determine which torus orbits.

    Let $\sigma$ be a codimension $k$ cone in the fan of $Y$.
    Then $V(\sigma)$ is a $k$-dimensional subscheme. 

    \begin{thm}
    \label{AlongOrbit}
    Let $w\in T_G^\vee$ be a point in the top dimensional cell of the Gr\"{o}bner fan. The multiplicity of
    $\init{w}{X}$ along $V(\sigma)$ is
    $$\sum_x m_x[T^\wedge:M_x+\sigma^\perp]$$
    where the sum is over
    all $x$ in $-\sigma^\circ\cap (-w+\Trop(X))$ and the 
    intersection multiplicities correspond to the
    intersection of $-w+\Trop(X)$ and $-\sigma$. 
    \end{thm}

    \begin{proof}
    
    We may refine $\Delta$ so that $X$ intersects torus orbits properly.  
    By the toric version of Chow's lemma, we may further refine $\Delta$ by so that $Y$ is smooth and projective.  Let $W$ be the complete intersection of $k$ ample hypersurfaces.  By applying the Kleiman-Bertini theorem on each torus orbit when choosing hypersurfaces, we may ensure that $W$ intersects torus orbits properly.  By ampleness,  $W\cap V(\sigma)\neq \emptyset$.
  
   $\Trop(W)$ is a union of cones of $\Delta$ of codimension at least $k$.   Let $d=\deg(W\cdot V(\sigma))$.   The multiplicity of the cone $-\sigma$ in $\Trop(W)$ is $d$.  By Lemma \ref{ConstantTransversality}, without changing $\Trop(W)$, we may replace $W$ by $\lambda\cdot W$ to ensure that $W$ intersects $t^w\cdot X$ in the interior.  If $Z$ is any components of $\init{w}{X}$ not supported on $V(\sigma)$, then $Z$ must intersect $V(\sigma)$ in a proper torus orbit.  Since $W$ intersects torus orbits properly, $W$ does not intersect $Z$ at any points of $V(\sigma)$.  

    Now $X\times_Y (t^{-w}\cdot W)$ is a zero-dimensional scheme supported on $T$.  Because specialization commutes with
    refined intersection product as in Theorem \ref{IntersectionNumber},
   \[\init{w}{X\cdot_Y (t^{-w}\cdot W)}=\init{w}{X}\cdot_{Y_0}
    \init{w}{t^{-w}\cdot W}=\init{w}{X}\cdot_{Y_0} W.\]
    We decompose the intersection
    product of $X$ and $t^{-w}\cdot W$ into contributions with different valuations as in the
    proof of Theorem \ref{IntersectionNumber}.
    Some contributions deform to give the intersection
    product of $\init{w}{X}$ and $W$ along the components of
    $\init{w}{X}$ supported on $V(\sigma)$.    By Lemma \ref{toruslimits}, $v(t^w x)\in \sigma^\circ$, if and only if $\init{w}{x}$ is
    a point in $\oh_\sigma$.  Therefore, the components of $X \cap (t^{-w}W)$ that deform to the intersection of $W$ with $V(\sigma)$ are the ones supported on $x$ with 
    \[v(x)\in (w-\Trop(X))\cap (w-\Trop(t^{-w}\cdot W))\cap \sigma^\circ=(w-\Trop(X))\cap\sigma^\circ.\]
    Each point $x$ counts with multiplicity $m_xd[T^\wedge:M_x+\sigma^\perp]$.  Since $\deg(W\cdot V(\sigma))=d$, we divide by $d$ to get the multiplicity of $\init{w}{X}$ along $V(\sigma)$.
    \end{proof}

    \bibliographystyle{plain}
    \bibliography{toolkit}

\vspace{+10 pt} \noindent

\end{document}